\magnification=1200


\hsize=14cm    
\vsize=20.5cm       
\parindent=0cm   \parskip=0pt     
\pageno=1 

\def\ind{\hskip 1cm\relax}


\hoffset=15mm    
\voffset=1cm    
 

\ifnum\mag=\magstep1
\hoffset=-0.5cm   
\voffset=-0.5cm   
\fi


\pretolerance=500 \tolerance=1000  \brokenpenalty=5000

\catcode`\@=11

\newcount\secno
\newcount\prmno
\newif\ifnotfound
\newif\iffound

\def\namedef#1{\expandafter\def\csname #1\endcsname}
\def\nameuse#1{\csname #1\endcsname}

\long\def\ifundefined#1#2#3{\expandafter\ifx\csname
  #1\endcsname\relax#2\else#3\fi}
\def\hwrite#1#2{{\let\the=0\edef\next{\write#1{#2}}\next}}

\toksdef\ta=0 \toksdef\tb=2
\long\def\leftappenditem#1\to#2{\ta={\\{#1}}\tb=\expandafter{#2}%
                                \edef#2{\the\ta\the\tb}}
\long\def\rightappenditem#1\to#2{\ta={\\{#1}}\tb=\expandafter{#2}%
                                \edef#2{\the\tb\the\ta}}

\def\lop#1\to#2{\expandafter\lopoff#1\lopoff#1#2}
\long\def\lopoff\\#1#2\lopoff#3#4{\def#4{#1}\def#3{#2}}

\def\ismember#1\of#2{\foundfalse{\let\given=#1%
    \def\\##1{\def\next{##1}%
    \ifx\next\given{\global\foundtrue}\fi}#2}}

\def\section#1{\medbreak
               \global\def\currenvir{section}
               \global\advance\secno by1\global\prmno=0
               {\bf \number\secno. {#1}}}

\def\subsection{\global\def\currenvir{subsection}
                \global\advance\prmno by1
                {\hskip-.5truemm(\number\secno.\number\prmno)}\hskip 5truemm}

\def\formule{\global\def\currenvir{formule}
                \global\advance\prmno by1
                {\hbox{\rm (\number\secno.\number\prmno)}}}

\def\proclaim#1{\global\advance\prmno by 1
                {\bf #1 \the\secno.\the\prmno.-- }}


\def\ex#1{\medbreak\global\def\currenvir{exemple}
																	\global\advance\prmno by 1
                {\bf #1 \the\secno.\the\prmno.}}

\def\qq#1{\medbreak\global\advance\prmno by 1
                { #1 \the\prmno)}}

\long\def\th#1 \enonce#2\endth{\global\def\currenvir{th}
   \medbreak\proclaim{#1}{\it #2}\medskip}

\long\def\comment#1\endcomment{}


\def\isinlabellist#1\of#2{\notfoundtrue%
   {\def\given{#1}%
    \def\\##1{\def\next{##1}%
    \lop\next\to\za\lop\next\to\zb%
    \ifx\za\given{\zb\global\notfoundfalse}\fi}#2}%
    \ifnotfound{\immediate\write16%
                 {Warning - [Page \the\pageno] {#1} No reference found}}%
                \fi}%
\def\ref#1{\ifx\labellist\empty{\immediate\write16
                 {Warning - No references found at all.}}
               \else{\isinlabellist{#1}\of\labellist}\fi}

\def\newlabel#1#2{\rightappenditem{\\{#1}\\{#2}}\to\labellist}
\def\labellist{}

\def\label#1{%
  \def\given{th}%
  \ifx\given\currenvir%
{\hwrite\lbl{\string\newlabel{#1}{\number\secno.\number\prmno}}}\fi%

\def\given{section}%
  \ifx\given\currenvir%
{\hwrite\lbl{\string\newlabel{#1}{\number\secno}}}\fi%

\def\given{subsection}%
  \ifx\given\currenvir%
{\hwrite\lbl{\string\newlabel{#1}{\number\secno.\number\prmno}}}\fi%

\def\given{formule}%
  \ifx\given\currenvir%
{\hwrite\lbl{\string\newlabel{#1}{\number\secno.\number\prmno}}}\fi%

\def\given{exemple}%
  \ifx\given\currenvir%
{\hwrite\lbl{\string\newlabel{#1}{\number\secno.\number\prmno}}}\fi%

\def\given{subsubsection}%
  \ifx\given\currenvir%
  {\hwrite\lbl{\string%
   
\newlabel{#1}{\number\secno.\number\subsecno.\number\subsubsecno}}}\fi
  \ignorespaces}

\newwrite\lbl

\def\openall{\openout\lbl=\jobname.lbl}
\def\closeall{\closeout\lbl}

\newread\testfile
\def\lookatfile#1{\openin\testfile=\jobname.#1
    \ifeof\testfile{\immediate\openout\nameuse{#1}\jobname.#1
                    \write\nameuse{#1}{}
                    \immediate\closeout\nameuse{#1}}\fi%
    \immediate\closein\testfile}%

\def\begin{\lookatfile{lbl}
           \input\jobname.lbl
           \openall}
\let\bye\end
\def\end{\closeall\bye}



\font\eightrm=cmr8         \font\eighti=cmmi8
\font\eightsy=cmsy8        \font\eightbf=cmbx8
\font\eighttt=cmtt8        \font\eightit=cmti8
\font\eightsl=cmsl8        \font\sixrm=cmr6
\font\sixi=cmmi6           \font\sixsy=cmsy6
\font\sixbf=cmbx6


\font\tengoth=eufm10       \font\tenbboard=msbm10
\font\eightgoth=eufm8      \font\eightbboard=msbm8
\font\sevengoth=eufm7      \font\sevenbboard=msbm7
\font\sixgoth=eufm6        \font\fivegoth=eufm5

\font\tencyr=wncyr10       
\font\eightcyr=wncyr8      
\font\sevencyr=wncyr7      
\font\sixcyr=wncyr6


\skewchar\eighti='177 \skewchar\sixi='177
\skewchar\eightsy='60 \skewchar\sixsy='60


\newfam\gothfam           \newfam\bboardfam
\newfam\cyrfam

\def\tenpoint{%
  \textfont0=\tenrm \scriptfont0=\sevenrm \scriptscriptfont0=\fiverm
  \def\rm{\fam\z@\tenrm}%
  \textfont1=\teni  \scriptfont1=\seveni  \scriptscriptfont1=\fivei
  \def\oldstyle{\fam\@ne\teni}\let\old=\oldstyle
  \textfont2=\tensy \scriptfont2=\sevensy \scriptscriptfont2=\fivesy
  \textfont\gothfam=\tengoth \scriptfont\gothfam=\sevengoth
  \scriptscriptfont\gothfam=\fivegoth
  \def\goth{\fam\gothfam\tengoth}%
  \textfont\bboardfam=\tenbboard \scriptfont\bboardfam=\sevenbboard
  \scriptscriptfont\bboardfam=\sevenbboard
  \def\bb{\fam\bboardfam\tenbboard}%
 \textfont\cyrfam=\tencyr \scriptfont\cyrfam=\sevencyr
  \scriptscriptfont\cyrfam=\sixcyr
  \def\cyr{\fam\cyrfam\tencyr}%
  \textfont\itfam=\tenit
  \def\it{\fam\itfam\tenit}%
  \textfont\slfam=\tensl
  \def\sl{\fam\slfam\tensl}%
  \textfont\bffam=\tenbf \scriptfont\bffam=\sevenbf
  \scriptscriptfont\bffam=\fivebf
  \def\bf{\fam\bffam\tenbf}%
  \textfont\ttfam=\tentt
  \def\tt{\fam\ttfam\tentt}%
  \abovedisplayskip=12pt plus 3pt minus 9pt
  \belowdisplayskip=\abovedisplayskip
  \abovedisplayshortskip=0pt plus 3pt
  \belowdisplayshortskip=4pt plus 3pt 
  \smallskipamount=3pt plus 1pt minus 1pt
  \medskipamount=6pt plus 2pt minus 2pt
  \bigskipamount=12pt plus 4pt minus 4pt
  \normalbaselineskip=12pt
  \setbox\strutbox=\hbox{\vrule height8.5pt depth3.5pt width0pt}%
  \let\bigf@nt=\tenrm       \let\smallf@nt=\sevenrm
  \normalbaselines\rm}

\def\eightpoint{%
  \textfont0=\eightrm \scriptfont0=\sixrm \scriptscriptfont0=\fiverm
  \def\rm{\fam\z@\eightrm}%
  \textfont1=\eighti  \scriptfont1=\sixi  \scriptscriptfont1=\fivei
  \def\oldstyle{\fam\@ne\eighti}\let\old=\oldstyle
  \textfont2=\eightsy \scriptfont2=\sixsy \scriptscriptfont2=\fivesy
  \textfont\gothfam=\eightgoth \scriptfont\gothfam=\sixgoth
  \scriptscriptfont\gothfam=\fivegoth
  \def\goth{\fam\gothfam\eightgoth}%
  \textfont\cyrfam=\eightcyr \scriptfont\cyrfam=\sixcyr
  \scriptscriptfont\cyrfam=\sixcyr
  \def\cyr{\fam\cyrfam\eightcyr}%
  \textfont\bboardfam=\eightbboard \scriptfont\bboardfam=\sevenbboard
  \scriptscriptfont\bboardfam=\sevenbboard
  \def\bb{\fam\bboardfam}%
  \textfont\itfam=\eightit
  \def\it{\fam\itfam\eightit}%
  \textfont\slfam=\eightsl
  \def\sl{\fam\slfam\eightsl}%
  \textfont\bffam=\eightbf \scriptfont\bffam=\sixbf
  \scriptscriptfont\bffam=\fivebf
  \def\bf{\fam\bffam\eightbf}%
  \textfont\ttfam=\eighttt
  \def\tt{\fam\ttfam\eighttt}%
  \abovedisplayskip=9pt plus 3pt minus 9pt
  \belowdisplayskip=\abovedisplayskip
  \abovedisplayshortskip=0pt plus 3pt
  \belowdisplayshortskip=3pt plus 3pt 
  \smallskipamount=2pt plus 1pt minus 1pt
  \medskipamount=4pt plus 2pt minus 1pt
  \bigskipamount=9pt plus 3pt minus 3pt
  \normalbaselineskip=9pt
  \setbox\strutbox=\hbox{\vrule height7pt depth2pt width0pt}%
  \let\bigf@nt=\eightrm     \let\smallf@nt=\sixrm
  \normalbaselines\rm}

\tenpoint


\def\pc#1{\bigf@nt#1\smallf@nt}         \def\pd#1 {{\pc#1} }


\catcode`\;=\active
\def;{\relax\ifhmode\ifdim\lastskip>\z@\unskip\fi
\kern\fontdimen2  -1.2 \fontdimen3 \string;}

\catcode`\:=\active
\def:{\relax\ifhmode\ifdim\lastskip>\z@\unskip\fi\penalty\@M\ \fi\string:}

\catcode`\!=\active
\def!{\relax\ifhmode\ifdim\lastskip>\z@
\unskip\fi\kern\fontdimen2  -1.1 \fontdimen3 \string!}

\catcode`\?=\active
\def?{\relax\ifhmode\ifdim\lastskip>\z@
\unskip\fi\kern\fontdimen2  -1.1 \fontdimen3 \string?}

\def\^#1{\if#1i{\accent"5E\i}\else{\accent"5E #1}\fi}
\def\"#1{\if#1i{\accent"7F\i}\else{\accent"7F #1}\fi}

\frenchspacing


\newtoks\auteurcourant      \auteurcourant={\hfil}
\newtoks\titrecourant       \titrecourant={\hfil}

\newtoks\hautpagetitre      \hautpagetitre={\hfil}
\newtoks\baspagetitre       \baspagetitre={\hfil}

\newtoks\hautpagegauche     
\hautpagegauche={\eightpoint\rlap{\folio}\hfil\the\auteurcourant\hfil}
\newtoks\hautpagedroite     
\hautpagedroite={\eightpoint\hfil\the\titrecourant\hfil\llap{\folio}}

\newtoks\baspagegauche      \baspagegauche={\hfil} 
\newtoks\baspagedroite      \baspagedroite={\hfil}

\newif\ifpagetitre          \pagetitretrue




\def\raggedbottom{\topskip 10pt plus 36pt\r@ggedbottomtrue}





\def\ctexte#1\endctexte{%
  \hbox{$\vcenter{\halign{\hfill##\hfill\crcr#1\crcr}}$}}


\long\def\ctitre#1\endctitre{%
    \ifdim\lastskip<24pt\vskip-\lastskip\bigbreak\bigbreak\fi
  		\vbox{\parindent=0pt\leftskip=0pt plus 1fill
          \rightskip=\leftskip
          \parfillskip=0pt\bf#1\par}
    \bigskip\nobreak}

\let\+=\tabalign

\def\signature#1\endsignature{\vskip 15mm minus 5mm\rightline{\vtop{#1}}}

\mathcode`A="7041 \mathcode`B="7042 \mathcode`C="7043 \mathcode`D="7044
\mathcode`E="7045 \mathcode`F="7046 \mathcode`G="7047 \mathcode`H="7048
\mathcode`I="7049 \mathcode`J="704A \mathcode`K="704B \mathcode`L="704C
\mathcode`M="704D \mathcode`N="704E \mathcode`O="704F \mathcode`P="7050
\mathcode`Q="7051 \mathcode`R="7052 \mathcode`S="7053 \mathcode`T="7054
\mathcode`U="7055 \mathcode`V="7056 \mathcode`W="7057 \mathcode`X="7058
\mathcode`Y="7059 \mathcode`Z="705A
 
\def\spacedmath#1{\def\packedmath##1${\bgroup\mathsurround=0pt ##1\egroup$}%
\mathsurround#1 \everymath={\packedmath}\everydisplay={\mathsurround=0pt }}
 
\def\nospacedmath{\mathsurround=0pt \everymath={}\everydisplay={} }


\def\decale#1{\smallbreak\hskip 28pt\llap{#1}\kern 5pt}
\def\decaledecale#1{\smallbreak\hskip 34pt\llap{#1}\kern 5pt}
\def\puce{\smallbreak\hskip 6pt{$\scriptstyle\bullet$}\kern 5pt}



\def\displaylinesno#1{\displ@y\halign{
\hbox to\displaywidth{$\@lign\hfil\displaystyle##\hfil$}&
\llap{$##$}\crcr#1\crcr}}


\def\ldisplaylinesno#1{\displ@y\halign{ 
\hbox to\displaywidth{$\@lign\hfil\displaystyle##\hfil$}&
\kern-\displaywidth\rlap{$##$}\tabskip\displaywidth\crcr#1\crcr}}


\def\eqalign#1{\null\,\vcenter{\openup\jot\m@th\ialign{
\strut\hfil$\displaystyle{##}$&$\displaystyle{{}##}$\hfil
&&\quad\strut\hfil$\displaystyle{##}$&$\displaystyle{{}##}$\hfil
\crcr#1\crcr}}\,}


\def\system#1{\left\{\null\,\vcenter{\openup1\jot\m@th
\ialign{\strut$##$&\hfil$##$&$##$\hfil&&
        \enskip$##$\enskip&\hfil$##$&$##$\hfil\crcr#1\crcr}}\right.}


\let\@ldmessage=\message

\def\message#1{{\def\pc{\string\pc\space}%
                \def\'{\string'}\def\`{\string`}%
                \def\^{\string^}\def\"{\string"}%
                \@ldmessage{#1}}}

\def\diagram#1{\def\normalbaselines{\baselineskip=0pt
\lineskip=10pt\lineskiplimit=1pt}   \matrix{#1}}



\def\up#1{\raise 1ex\hbox{\smallf@nt#1}}


\def\cf{{\it cf}} 

\def\qed{\raise -2pt\hbox{\vrule\vbox to 10pt{\hrule width 4pt
                 \vfill\hrule}\vrule}}

\def\cqfd{\unskip\penalty 500\quad\vrule height 4pt depth 0pt width 4pt\medbreak}

\def\virg{\raise .4ex\hbox{,}}   


\def\build#1_#2^#3{\mathrel{
\mathop{\kern 0pt#1}\limits_{#2}^{#3}}}


\def\boxit#1#2{%
\setbox1=\hbox{\kern#1{#2}\kern#1}%
\dimen1=\ht1 \advance\dimen1 by #1 \dimen2=\dp1 \advance\dimen2 by #1 
\setbox1=\hbox{\vrule height\dimen1 depth\dimen2\box1\vrule}%
\setbox1=\vbox{\hrule\box1\hrule}%
\advance\dimen1 by .4pt \ht1=\dimen1 
\advance\dimen2 by .4pt \dp1=\dimen2  \box1\relax}

\def\date{\the\day\ \ifcase\month\or janvier\or f\'evrier\or mars\or
avril\or mai\or juin\or juillet\or ao\^ut\or septembre\or octobre\or
novembre\or d\'ecembre\fi \ {\old \the\year}}

\def\dateam{\ifcase\month\or January\or February\or March\or
April\or May\or June\or July\or August\or September\or October\or
November\or December\fi \ \the\day ,\ \the\year}

\def\moins{\mathop{\hbox{\vrule height 3pt depth -2pt
width 5pt}}}
\def\crog{{\vrule height 2.57mm depth 0.85mm width 0.3mm}\kern -0.36mm
[}

\def\crod{]\kern -0.4mm{\vrule height 2.57mm depth 0.85mm
width 0.3 mm}}
\def\moins{\mathop{\hbox{\vrule height 3pt depth -2pt
width 5pt}}}

\def\rond{\kern 1pt{\scriptstyle\circ}\kern 1pt}

\def\diagram#1{\def\normalbaselines{\baselineskip=0pt
\lineskip=10pt\lineskiplimit=1pt}   \matrix{#1}}

\def\hfl#1#2{\nospacedmath\smash{\mathop{\hbox to
12mm{\rightarrowfill}}\limits^{\scriptstyle#1}_{\scriptstyle#2}}}

\def\ghfl#1#2{\nospacedmath\smash{\mathop{\hbox to
25mm{\rightarrowfill}}\limits^{\scriptstyle#1}_{\scriptstyle#2}}}

\def\phfl#1#2{\nospacedmath\smash{\mathop{\hbox to
8mm{\rightarrowfill}}\limits^{\scriptstyle#1}_{\scriptstyle#2}}}

\def\pvfl#1#2{\llap{$\scriptstyle#1$}\left\downarrow\vbox to
3.5mm{}\right.\rlap{$\scriptstyle#2$}}

\def\pa{\S\kern.15em}

\def\Z{{\bf Z}}

\def\P{{\bf P}}

\def\Q{{\bf Q}}
\def\C{{\bf C}}

\def\cad{c'est-\`a-dire}

\def\Card{\mathop{\rm Card}\nolimits}

\def\cf{{\it cf.\/}}

\def\Hom{\mathop{\rm Hom}\nolimits}

\def\Id{\hbox{\rm Id}}

\def\Im{\mathop{\rm Im}\nolimits}

\def\isom{\simeq}

\def\Ker{\mathop{\rm Ker}\nolimits}

\def\loc{{\it loc.cit.\/}}
\def\long{\mathop{\rm long}\nolimits}
\def\lra{\longrightarrow}
\def\llra{\nospacedmath\hbox to 10mm{\rightarrowfill}}
\def\lllra{\nospacedmath\hbox to 15mm{\rightarrowfill}}

\def\Pic{\mathop{\rm Pic}\nolimits}

\def\ra{\rightarrow}

\def\rang{\mathop{\rm rang}\nolimits}

\def\rang{\mathop{\rm rg}\nolimits}

\def\ssi{si et seulement si}
\def\theo{th\'eor\`eme}

\def\tx{\kern -1.5pt -}
\def\vide{\varnothing}

\def\cc#1{\hfill\kern .7em#1\kern .7em\hfill}

\def\ndeg{non d\'eg\'en\'er\'e}

\def\note#1#2{\footnote{\parindent
0.4cm$^#1$}{\vtop{\eightpoint\baselineskip12pt\hsize15.5truecm
\noindent #2}}\parindent 0cm}

\def\og{\leavevmode\raise.3ex\hbox{$\scriptscriptstyle\langle\!\langle$}}
\def\fg{\leavevmode\raise.3ex\hbox{$\scriptscriptstyle\,\rangle\!\rangle$}}

\def\a{{\alpha}}
\def\b{{\beta}}
\def\d{{\delta}}

\def\epsilon{{\varepsilon}}
\def\eps{{\varepsilon}}

\def\l{{\lambda}}

\def\s{{\sigma}}

\def\L{\Lambda}
\def\cS{{\cal S}}

\def\cF{{\cal F}}

\def\cH{{\cal H}}

\def\cK{{\cal K}}

\def\cP{{\cal P}}

\def\cO{{\cal O}}

\catcode`\@=12

\showboxbreadth=-1  \showboxdepth=-1


\spacedmath{1.5pt}\parskip=1mm 
\baselineskip=14pt
\overfullrule=0pt
\input amssym.def
\input amssym
 
\def\tilde{\widetilde}

\def\Gr{\mathop{\rm Gr}\nolimits}
\def\tilde{\widetilde} 
\def\eps{\epsilon}
\def\vide{\varnothing}
\def\F{{\bf F}}

\def\lda#1#2{\llap{$\scriptstyle#1$}\left\downarrow\vbox to
3mm{}\right.\rlap{$\scriptstyle#2$}}

\null\bigskip
\centerline {\bf LIEUX DE D\'EG\'EN\'ERESCENCE}

\centerline{\bf O. Debarre}

\vskip 1cm 

\begin

\ind Soient $X$ une vari\'et\'e complexe projective {\it lisse} connexe
et
$Y$ le lieu des z\'eros d'une section d'un fibr\'e en droites ample sur
$X$. Le
\theo\ de Lefschetz \'enonce que la restriction $H^p(X,\Z)\to H^p(Y,\Z)$
est bijective pour $p<\dim X-1$, injective pour $p=\dim X-1$. Cela
entra\^ine   le
\theo\ de Bertini: $Y$ est connexe si sa dimension est au moins $1$.
Dans le m\^eme ordre d'id\'ees, Grothendieck a montr\'e que les groupes
de Picard de
$X$ et de
$Y$ sont isomorphes (quelles que soient les singularit\'es de $Y$!) si
$\dim Y\ge 3$.
 
\ind Etant donn\'es des 
fibr\'es vectoriels
$E$ et
$F$ sur $X$ de rangs respectifs $e$ et $f$,  et un morphisme $u:E\ra F$,
on consid\`ere les lieux de d\'eg\'en\'erescence
$D_r=\{ x\in X\mid \rang (u_x)\le r\}$.
 Fulton et Lazarsfeld ont d\'emontr\'e dans [FL] l'analogue du
\theo\ de Bertini : si $\Hom(E,F)$ est
ample, $D_r$ est connexe si sa dimension attendue $\d(r)=
\dim(X)-(f-r)(e-r)$ est au moins $1$. Nous poursuivons leurs
m\'ethodes pour obtenir des extensions  des \theo s de Lefschetz et
Grothendieck mentionn\'es plus haut. Il y a plusieurs cas de
figure:

\ind $\bullet$ si $D_{r-1}$ est vide, on peut compl\`etement d\'ecrire
(\cf\ (\ref{cons})) la cohomologie enti\`ere de
$D_r$ jusqu'en degr\'e $\d(r)-1$. Si
$D_r$ est normal, que $0<r<\min \{e,f\}$ et
$\d(r)\ge 3$, le groupe de Picard de
$D_r$ est isomorphe \`a $\Pic (X)\oplus \Z$.

\ind $\bullet$ Si au contraire suffisament des lieux $D_s$, pour $s\le
r$, sont non vides, la
restriction
$H^p(X,\Z)\to H^p(D_r,\Z)$ est  un isomorphisme pour $p$ assez petit
(\cf\ th.~\ref{lef} pour un \'enonc\'e pr\'ecis). Par exemple, si
$D_r$ est normal, que  
$D_{r-1}$ n'est pas vide ou que $r=0$, et que $\d(r)\ge 3$, les
groupes de Picard de
$X$ et de $D_r$ sont isomorphes (cor.~\ref{cor34}).

\ind Dans [E], Ein montre ces r\'esultats sur le groupe de
Picard dans le cas o\`u
$E$ est trivial,  en supposant seulement
$\d(r)=2$ mais aussi $F$ \og suffisamment ample\fg\  et $u$ g\'en\'eral
(c'est une extension du \theo\ de Noether-Lefschetz sur le groupe de
Picard d'une surface g\'en\'erale dans $\P^3$). 

\ind D'autre part, pour
$r=\min\{e,f\}-1$ et $D_{r-1}$ vide, toujours sous les
hypoth\`eses
$F$
\og suffisamment ample\fg\  et
$u$ g\'en\'eral, Spandaw d\'etermine dans sa th\`ese les classes
{\it alg\'ebriques} de
$H^{\d(r)}(D_r,\Z)$. 
 
\ind Nous nous int\'eressons ensuite au cas d'un morphisme $u:E\to
E^*\otimes L$ {\it antisym\'etrique}. Tu a d\'emontr\'e que si
$\wedge^2E^*\otimes L$ est ample,  le lieu de d\'eg\'en\'erescence
$A_r=\{ x\in X\mid \rang (u_x)\le 2r\}$ est connexe si sa dimension
attendue  $\a(r)=\dim(X)-{e-2r\choose 2}$ est au moins $1$. Nous 
obtenons des extensions  des \theo s de 
Grothendieck et Lefschetz dans ce cadre:
la restriction $\Pic (X)\to \Pic (A_r)$ est bijective si
$A_r$ est normale et $\a(r)\ge 3$, et 

\ind $\bullet$ si $A_{r-1}$ est vide, on peut d\'ecrire (th.~\ref{ortho}) la
cohomologie enti\`ere de
$A_r$ jusqu'en degr\'e $\a(r)-1$;

\ind $\bullet$ si au contraire suffisament des lieux $A_s$, pour $s\le
r$, sont non vides (\cf\ th.~\ref{lefalt} pour un \'enonc\'e pr\'ecis),
la restriction
$H^p(X,\Z)\to H^p(A_r,\Z)$ est un isomorphisme pour $p$ assez petit.

\ind Nous terminons par l'\'etude du cas des fibr\'es {\it
orthogonaux} (dont le cas des morphismes antisym\'etriques est un cas
particulier): on se donne  un fibr\'e vectoriel
$V$ de rang 
pair sur
$X$ muni d'une forme quadratique
\ndeg e \`a valeurs dans un fibr\'e en droites $L$ et des 
sous-fibr\'es $E$ et $F$ totalement isotropes maximaux de $V$. On montre
un \theo\ de Bertini pour les lieux  de d\'eg\'en\'erescence
$$ O^r=\{\ x\in X\mid \dim(E_x\cap F_x)\ge  r\ \ {\rm et}\ \
\dim(E_x\cap F_x)\equiv r\pmod{2}\ \}\ ;$$ si $E^*\otimes F^*\otimes L$
est ample,
$ O^r$ est connexe si sa dimension attendue
$\a(r)=\dim(X)-{r\choose 2}$ est au moins $1$. Les r\'esultats de type
Lefschetz obtenus dans le cas antisym\'etrique devraient subsister, mais nous
n'en obtenons qu'une maigre confirmation (prop.~\ref{OO}). 

\ind Dans cet article, tous les sch\'emas
sont de type fini sur le corps des nombres complexes. On d\'esigne
par $\F$ un corps fini ou \'egal \`a $\Q$.

\ind Je remercie R. Laterveer et W. Fulton pour leur aide pour les
prop.~\ref{rob} et \ref{prop} respectivement.
\vskip1cm
 \centerline{\bf I. Lieux de d\'eg\'en\'erescence}
\bigskip

\section{Le r\'esultat de Fulton et Lazarsfeld}

\ind Soient $X$ une vari\'et\'e complexe projective irr\'eductible  et
$E$ et
$F$ des fibr\'es vectoriels sur $X$ de rangs respectifs $e$ et $f$. Soit 
$u:E\ra F$ un morphisme; on note 
$$D_r=\{ x\in X\mid \rang (u_x)\le r\}
$$ et on pose
$\d(r)=
\dim(X)-(f-r)(e-r)$. Par la suite, nous supposerons toujours $ e\le f$ (ce
que l'on peut toujours faire quitte \`a remplacer $u$ par son dual).
 
\subsection Soient $\pi:G=G(e-r,E)\to X$ le fibr\'e en
grassmanniennes et
$S$ le fibr\'e tautologique de rang $e-r$ sur $G$. Soit $Y$ le lieu
des z\'eros de la compos\'ee\label{cons}
$$ S\hookrightarrow \pi^*E\buildrel{\pi^*u}\over{\lra}\pi^*F\ .$$
\ind Le morphisme $\pi$ induit par restriction un morphisme
$\pi':Y\ra D_r$ propre surjectif, birationnel au-dessus de $D_r\moins
D_{r-1} $, de fibre $G(e-r,e-l)$ au-dessus de $D_l\moins D_{l-1} $. 
 Fulton et Lazarsfeld montrent  que {\it si  
$\Hom(E,F)$ est ample},
$H^q(G\moins Y,\F) $ s'annule pour
$q\ge\dim(X)+(f+r)(e-r)$. Par dualit\'e de Lefschetz, on en d\'eduit,
 si $X$ est lisse, 
$H^p(G,Y;\F)=0$ pour
$p\le\d(r)$. 

 \ind La dualit\'e de Lefschetz n'est valable que lorsque
$G\moins Y$ est lisse. Elle est remplac\'ee dans le cas g\'en\'eral par une
suite spectrale
$$E_2^{pq}=H^p(G\moins Y,\cH_{-q}(G,\F))\Rightarrow H_{-p-q}(G,Y;\F)\ ,$$
 o\`u $\cH_q(G,\F)$ est le faisceau de fibre $H_q(G,G\moins\{x\};\F)$ en
un point $x$ de $G$ ([H1], p.~548). Lorsque $G\moins Y$ est localement
intersection compl\`ete,   le
support de
$\cH_{\dim(G)+i}
(G,\F)$ est de dimension au plus $i$, pour tout $i\in\Z$ ([H1], lemma 4,
p.~550). Or la d\'emonstration de Fulton et Lazarsfeld montre que pour tout
ferm\'e $Z$ de
$X$, la dimension cohomologique de $ \pi^{-1}(Z)\moins Y$ est au plus
$\dim(\pi^{-1}(Z))+f(e-r)-1$. On en d\'eduit\note{1}{Les faisceaux $\cH_q$ ne
sont localement constant que sur chaque strate d'une stratification de
Whitney de $G$, et il faut en fait raisonner strate par strate comme dans 
[H2], Lemma~3, p.~134.}
$$E^2_{pq}=0\qquad{\rm pour}\quad p>-q-\dim(G)+f(e-r)-1\ ,$$
d'o\`u de nouveau 
$$H_p(G,Y;\F)=0 \qquad{\rm pour}\quad p\le\dim(G)-f(e-r)=\d(r)\
,\leqno{\formule}\hbox{\label{FLL}}$$ {\it sous l'hypoth\`ese que
$X\moins D_0$ est 
 localement
intersection compl\`ete}\note{2}{On a une autre suite
spectrale 
 $$E_2^{pq}=H^p(G,Y;\cH_{-q}(G,\F))\Rightarrow H_{-p-q}(G\moins Y,\F)\ ,$$
qui permet de montrer que l'on a une inclusion $H^1(G,Y;\Q)\hookrightarrow
H_{2\dim(G)-1}(G\moins Y,\Q)$ sous la seule hypoth\`ese que $G$ est normale
([FL], lemma~1.3). Il en r\'esulte que $D_r$ est connexe d\`es que $\d(r)>0$,
sans hypoth\`ese sur les singularit\'es de $X$. Cela laisse \`a penser que
l'on doit pouvoir am\'eliorer l'hypoth\`ese \og $X $   
 localement
intersection compl\`ete\fg .}.

\medskip
\subsection Notons $\iota_r$ l'inclusion $D_r\hookrightarrow X$ et $c$
l'inverse dans
$H^{\bullet}(D_r,\Z )$ de la classe de Chern totale du fibr\'e vectoriel
$\Ker(u)\vert_{D_r}$. Lorsque $D_{r-1}$ est vide, la discussion ci-dessus
entra\^ine que l'application\label{FL} 
$$\matrix{\displaystyle\bigoplus_{\l=(\l_1,\ldots,
\l_{e-r})\atop r\ge\l_1\ge\cdots\ge \l_{e-r}\ge 0}\hskip-5mm
H^{p-2|\l|}(X,\F )&\lra&H^p(D_r,\F )\cr \sum_\l
\a_\l&\longmapsto&\sum_\l 
\Delta_\l(c)\cdot \iota_r^*\a_\l\cr}$$
 est injective pour $p\le\d(r)$,
bijective pour $p<\d(r)$. On a employ\'e les notations standard
 $$|\l|=\l_1+\cdots+\l_{e-r}\hskip 9mm,\hskip 9mm
\Delta_\l(c)=\det(c_{\l_i+j-i})_{1\le i,j\le {e-r}} \ .$$

\ex{Exemples.} 1) {\bf Vari\'et\'es de Segre.} Soient $V$ et $W$ des
espaces vectoriels. On consid\`ere dans $\P\Hom(V,W)$ le lieu $D_1$
associ\'e au   morphisme tautologique
$u:V\otimes\cO \to W\otimes\cO(1)$. C'est l'image du
plongement de Segre $\P=\P  V^* \times\P  W \hookrightarrow\P \Hom(V,W)
$, ou encore le lieu des homomorphismes de rang $1$, et il a  la
dimension attendue $\d(1)$. Posons $h_1=c_1(p_1^*\cO_{\P V^*}(1))$ et
$h_2=c_1(p_2^*\cO_{\P W}(1))$; le fibr\'e $K=\Ker(u)$ sur $\P$
  s'ins\`ere dans une suite exacte\label{segre}
$$0\to K\to \cO_\P\otimes  V\to
p_2^*\cO_{\P W}(1)\to 0\ ,$$ de sorte que $1/c(K)=1 +h_2$; puisque $D_0$
est vide, (\ref{FL}) entra\^ine que pour $p<\d(1)$, le groupe
$H^p(\P,\Z)$ est nul pour
$p$ impair et a pour base  
$(h_2^s(h_1+h_2)^{p/2- s})_{0\le s\le p/2}$ pour $ p$ pair, ou encore  
$(  h_1^sh_2^{p/2- s})_{0\le s\le p/2}$.
 
2) Soit $C$ une
courbe projective lisse de genre
$g$. La vari\'et\'e 
$W^s_d(C)$ peut s'interpr\'eter comme le lieu $D_r$ pour un morphisme
$u:E\ra F$ de fibr\'es vectoriels sur la jacobienne de $C$, avec $r=e-s-1$
et $f =e+g-1-d$. Supposons
$W^1_d(C)$ vide, de sorte que
$W_d(C)$ est isomorphe au produit sym\'etrique $C_d$; on d\'eduit de
(\ref{FL}) un isomorphisme
$$H^p(C_d,\Z )\isom \bigoplus_{0\le 2j\le p}
x^j\cdot\iota^*H^{p-2j}(JC,\Z ) 
$$
pour $p<d$, avec $x=[C_{d-1}]$. C'est un cas particulier des formules de
MacDonald ([M]).
De la m\^eme fa\c con, le \theo\ permet de calculer $H^p(W^s_d(C),\Z)$
pour $p< g-(s+1)(g-d+s)$ lorsque $C$ n'a pas de $g^{s+1}_d$.

\bigskip

\ind  Nous aurons besoin d'une g\'en\'eralisation (bas\'ee
sur les id\'ees de [S]) du r\'esultat de connexit\'e de
Fulton et Lazarsfeld, qui fait intervenir la notion de $d$-connexit\'e. 
Rappelons qu'un sch\'ema $X$ est dit $d$\tx connexe s'il est de dimension
$>d$ et si, pour tout sous-sch\'ema ferm\'e $Z$ de $X$ de dimension $<d$, le
sch\'ema
$X\moins Z$ est connexe. Les propri\'et\'es  suivantes sont
classiques:

\ind 1) un sch\'ema est $(-1)$\tx connexe si et seulement s'il est
non vide; il est $0$\tx connexe si et seulement s'il est  connexe.

\ind 2) Un sch\'ema irr\'eductible de dimension $d$ est
$(d-1)$\tx connexe. Toute composante irr\'eductible d'un sch\'ema $d$\tx
connexe est de dimension $>d$.

\ind 3) Si $X$ est r\'eunion de sous-sch\'emas ferm\'es $d$\tx
connexes $X_1,\ldots,X_m$, il est $d$\tx connexe \ssi , pour tous $i$ et
$j$, il existe des indices
$i_0,i_1,\ldots,i_m$ avec $i_0=i$ et $i_m=j$ tels que
$\dim(X_{i_\nu}\cap X_{i_{\nu+1}})\ge d$ pour tout $\nu=0,\ldots,m-1$.

\th\label{FL2}
Proposition
\enonce
Soient $X$ un sch\'ema projectif $d$\tx connexe, et $E$ et
$F$ des fibr\'es vectoriels sur $X$ de rangs respectifs $e$ et $f$, avec
$\Hom(E,F)$ ample. Soit $u:E\ra F$ un 
 morphisme  partout de rang $\le k$; pour tout $r\le k$, le lieu 
$D_r$ est  $(d-(f-r)(e-r)+(e-k)(f-k))$\tx connexe.
 En particulier, si $X$ est irr\'eductible de dimension $>
(f-r)(e-r)-(e-k)(f-k)$, le lieu $D_r$ est connexe.
 \endth

{\it D\'emonstration.} Il suffit de traiter le cas $k=r+1$; posons 
$$d'=d-(f-r)(e-r)+(e-k)(f-k)=d-e-f+2r+1\ .$$
\ind Notons
$X_1,\ldots,X_m$ les composantes irr\'eductibles de $X$; elles sont
toutes de dimension $>d$ par 2). Par [S] Lemma 4.1.3 (qui se
d\'emontre aussi \`a partir de [ACGH], prop. (1.3), p. 307, en
prenant des sections hyperplanes),
chaque intersection
$X_j\cap D_r$ est
$d'$\tx connexe. Pour tous $i$ et $j$, il existe par 3) des indices
$i_0=i,i_1,\ldots,i_m=j$ tels que $\dim(X_{i_\nu}\cap X_{i_{\nu+1}})\ge d$
pour tout $\nu=0,\ldots,m-1$. On a par \loc
$$\dim(X_{i_\nu}\cap X_{i_{\nu+1}}\cap D_r)\ge d'\ ,$$
pour tout $\nu$, de sorte que $D_r=\cup_j (X_j\cap D_r)$ est $d'$\tx
connexe par 3).\cqfd

\section {Un \theo\ de Lefschetz}
 
\subsection Restons dans la situation du \S 1, dont nous gardons les
notations. Comme  remarqu\'e dans [FL], l'application
$H^p(\iota_r,\F)$ est injective pour
$p\le\d(r)$. On posera
$\eps(0)=1$,
$\eps(1)=2$ et, pour tout entier
$k$ strictement positif, $\eps(2k)=0$ et $\eps(2k+1)=1$.\label{rem} 
  
\th\label{lef}
Th\'eor\`eme
\enonce
Soit $X$ une vari\'et\'e projective irr\'eductible localement
intersection compl\`ete. Soient 
$E$ et
$F$ des fibr\'es vectoriels sur $X$, avec
$\Hom(E,F)$ {\rm ample}, et $u:E\to F$ un morphisme.
Supposons $[{m\over 2}]\le r$ et $\d(r-[{m\over 2}])\ge  \eps(m)$;
l'application
$H^p(\iota_r,\Z)$ est bijective pour $p\le
m$.\endth

{\it D\'emonstration.} Elle consiste \`a comparer les suites spectrales de
Leray pour l'application $\pi:G\to X$ et sa restriction
$\pi':Y\to D_r$. Les fibres de ces deux applications
\'etant des grassmanniennes, les faisceaux $R^q\pi_*\Z$ et $R^q\pi'_*\Z$ sont nuls
pour $q$ impair, de sorte que
${}^\pi\!E_2^{p,q}={}^\pi\!E_3^{p,q}$ et
${}^{\pi'}\!E_2^{p,q}={}^{\pi'}\!E_3^{p,q}$. D'autre part, la suite
spectrale
$${}^\pi\!E_2^{pq}=H^p(X,R^q\pi_*\Z)\ \ \Rightarrow\ \ H^{p+q}(G,\Z)$$
d\'eg\'en\`ere. Soit $x$ un point de $X$, on a
$$(R^{2q}\pi_*\Z)_x\isom H^{2q}(G(e-r,e),\Z)\isom 
\bigoplus_{r\ge\l_1\ge\cdots\ge
\l_{e-r}\ge 0\atop \l_1+\cdots+
\l_{e-r}=q}\Z
$$
et, si $x\in D_l\moins D_{l-1} $, 
$$(R^{2q}\pi'_*\Z)_x\isom H^{2q}(G(e-r,e-l),\Z)\isom 
\bigoplus_{r-l\ge\l_1\ge\cdots\ge
\l_{e-r}\ge 0\atop \l_1+\cdots+
\l_{e-r}=q}\Z
\ .$$
\ind En particulier, la restriction $R^{2q}\pi_*\Z\to R^{2q}\pi'_*\Z$ est
surjective pour tout $q$.  Son noyau $\cK_{2q}$ est nul sur $D_{r-q}$
et constant sur chaque $D_l\moins D_{l-1}$. Plus pr\'ecis\'ement, il
existe une filtration 
$$0=\cF_0\i\cF_1\i\cdots\i\cF_q\i\cF_{q+1}=\cK_{2q}\ ,$$ v\'erifiant
$\cF_{i+1}/\cF_i\isom \Z^{r_i}_{X\moins D_{r-q+i}}$ pour $0\le i\le
q$, o\`u les $r_i$ sont des entiers.

\ind  L'assertion \`a d\'emontrer d\'ecoule de (\ref{FLL}) pour $m=0$;
supposons-la vraie pour tout entier $<m$. On supposera aussi $r<e$; on
v\'erifie que dans ce cas, on a pour tous entiers $t\ge s\ge 0$ les
in\'egalit\'es
$$\d(t)\ge \d(t-s)+s(s+2)\leqno{\formule}\hbox{\label{formule1}}$$
\vskip-8mm$$\eps(t)+[{t\over 2}]([{t\over 2}]+2)>t\
.\leqno{\formule}\hbox{\label{formule2}}$$
{\bf Premier pas.} {\it Supposons $0\le p<m$ et $q+[{p\over 2}] \le
r$. On a} 
$$\eqalign{H^0(\cK_{2q})=0\qquad &{\rm si} \quad\d(r-q)\ge 0{\rm\ ;}\cr
H^{p+1}(\cK_{2q})=0\qquad &{\rm si} \quad \d(r-q-[{p\over 2}])\ge 
\eps(p)\ .\cr}$$
\ind
Consid\'erons la suite de cohomologie associ\'ee
\`a la suite exacte  
$$0\to\Z_{X\moins D_{r-q+i}}\to \Z_X\to \Z_{D_{r-q+i}}\to 0\ .$$
\ind Le \theo\ de Fulton et Lazarsfeld (\cf\ (\ref{FLL})) entra\^ine que
$H^0(
\Z_{X\moins D_{r-q+i}})$ est nul lorsque $\d(r-q+i)\ge 0$.  Comme la
fonction $\d$ est croissante, on en d\'eduit $H^0(\cK_{2q})=0$ lorsque
$\d(r-q)\ge 0$. 

\ind D'autre part, on a par
(\ref{rem}) des suites exactes
$$0\to H^p(\Z_X)\to H^p(\Z_{D_{r-q+i}})\to H^{p+1}(\Z_{X\moins
D_{r-q+i}})\to 0$$ lorsque $p<\d(r-q+i)$; l'hypoth\`ese de r\'ecurrence
entra\^ine alors
$H^{p+1}(\Z_{X\moins D_{r-q+i}})=0$ si l'on a de plus
$\d(r-q+i-[{p\over 2}])\ge \eps(p)
$ et $p<m$. Comme la fonction $\d$ est croissante, on a donc
$H^{p+1}(\cK_{2q})=0$ lorsque $\d(r-q-[{p\over 2}])\ge  \eps(p)$ et
$p<m$, puisque l'on a alors, par (\ref{formule1}) et (\ref{formule2}),
$$ \d(r-q )\ge \d(t-q-[{p\over 2}])+[{p\over 2}]([{p\over 2}]+2)
\ge \eps(p)+[{p\over 2}]([{p\over 2}]+2)>p\ .$$
\ind Ceci montre le premier pas.
\medskip

{\bf Deuxi\`eme pas.} {\it Supposons $[{m\over 2}]\le r$ et
$\d(r-[{m\over 2}])\ge 
\eps(m)
$.  L'application naturelle\break
$\iota_\infty^{p,q}:{}^\pi E_\infty^{p,q}\to
{}^{\pi'}\!E_\infty^{p,q}$ est injective pour  $p<m$ et
$p+q\le m$. D'autre part,\break  
${}^{\pi'}\!E_\infty^{m,0}\isom {}^{\pi'}\!E_2^{m,0}\isom
H^m(D_r,\Z)$.}
 
\ind  Supposons $p+2q\le m$. L'hypoth\`ese
$\d(r-[{m\over 2}])\ge  \eps(m) $ entra\^ine 
$$\d(r-q)\ge\d(r-[{m\over 2}])\ge  \eps(m)\ge 0$$ et, si $p+2q<m$,
$$\d(r-q-[{p\over 2}])\ge  \eps(p) \ ;$$
en effet,  si $q+[{p\over 2}]<[{m\over 2}]$, cela d\'ecoule de 
(\ref{formule1}), et si $q+[{p\over 2}]=[{m\over 2}]$, on a
$\eps(p)\le\eps(m)$.
 On consid\`ere le diagramme commutatif:
$$\diagram{H^p(\cK_{2q})\cr
\pvfl{}{}\cr
H^p(X,R^{2q}\pi_*\Z)&\hfl{{}^\pi\!d_3^{p,2q}
}{} &H^{p+3} (X,R^{2q-2}
\pi_*\Z)\cr
\pvfl{}{\iota_3^{p,2q}}&&\pvfl{}{\iota_3^{p+3,2q-2}} \cr
H^p(D_r,R^{2q}\pi'_*\Z)&\hfl{{}^{\pi'}\!\!d_3^{p,2q}
}{} &\ H^{p+3}(D_r,R^{2q-2}
\pi'_*\Z)
\cr\pvfl{}{}\cr
H^{p+1}(\cK_{2q})\ ; \cr}$$
il r\'esulte du premier pas que l'application
$\iota_3^{p,2q}$ est injective, et bijective si $p+2q<m$. 

\ind  Supposons maintenant $p+q\le m$; comme  ${}^\pi\! d_3^{p,q}=0$,
on a  
${}^{\pi'}\!\!d_3^{p,q}=0$  pour $p+q<m$. En particulier, 
${}^{\pi'}\!E_4^{p,q}$ est un sous-groupe de $ 
{}^{\pi'}\!E_3^{p,q}$ qui lui est \'egal pour $p=m$
et, le m\^eme raisonnement
s'appliquant
\`a chaque cran de la suite  spectrale,   
${}^{\pi'}\!E_\infty^{p,q}$ est un sous-groupe de $ 
{}^{\pi'}\!E_3^{p,q}$ qui lui est \'egal pour $p=m$, de sorte que
$\iota_\infty^{p,q}$ est injective.
\medskip

{\bf Conclusion.} Supposons $p+q\le m$ et
$\d(r-[{m\over 2}])\ge\eps(m)$; on en d\'eduit $\d(r)> m$ par
(\ref{formule1}) et (\ref{formule2}), de sorte que la restriction 
$\rho:H^{p+q}(G,\Z)\to H^{p+q}(Y,\Z)$ est bijective par
(\ref{FLL}). Comme $
\iota_\infty^{p,q}$ est injective (deuxi\`eme pas), il en
r\'esulte que
$\Gr\rho$ est bijective, ainsi donc que $
\iota_\infty^{m,0}$, qui, par le deuxi\`eme pas, n'est autre que
$H^m(\iota_r,\Z)$.\cqfd

\medskip
\ex{Remarque}  On a  
$R^2\pi'_*\Z\isom\Z_{D_{r-1}}$, d'o\`u un diagramme
commutatif \`a lignes  exactes, o\`u tous les groupes de
cohomologie sont \`a coefficients dans $\Z$\label{rem23}
$$\diagram{
0&\hskip-3mm\to\hskip-3mm&H^2(X)&\phfl{H^2(\pi)}{}&H^2(G)&\hskip-3mm
\to\hskip-3mm& H^0(X)&\hskip-3mm\buildrel{0}\over{\lra}\hskip-3mm
&H^3(X)&\phfl{H^3(\pi)}{}
&H^3(G)&\hskip-3mm\to\hskip-3mm& H^1(X)\cr
\noalign{\vskip 3mm}
&&\pvfl{}{\hskip-1mm
H^2(\iota_r)}&&\pvfl{}{}&&\pvfl{}{\hskip-1mm H^0(\iota_{r-1})}&&
\pvfl{}{\hskip-1mm H^3(\iota_r)}&&\pvfl{}{}&&\pvfl{}{\hskip-1mm
H^1(\iota_{r-1})}\cr
\noalign{\vskip 3mm}
 0&\hskip-3mm\to\hskip-3mm&
H^2(D_r)&\phfl{H^2(\pi')}{}&H^2(Y)&\hskip-3mm\to\hskip-3mm&
H^0(D_{r-1})& \phfl{{}^{\pi'}\!\!d_3^{02}}{} 
&H^3(D_r)&\phfl{H^3(\pi')}{}
&H^3(Y)&\hskip-3mm\to\hskip-3mm&  H^1(D_{r-1})\cr} 
$$ 
\ind Supposons $0<r<e$ et $\d(r)\ge 3$. Si  $D_{r-1}$ n'est pas vide,
on en d\'eduit le diagramme 
$$\diagram{
&&0&\to&H^3(X)&\buildrel{\pi^*}\over{\lra}&H^3(G)
&\to &H^1(X)\ \cr
\noalign{\vskip 2mm}
&&\lda{}{}&&\lda{}{\hskip-1mm
H^3(\iota_r)}&&\lda{}{}&&\lda{}{\hskip-1mm
H^1(\iota_{r-1})}\cr 0&\to&\Z^{c-1}&\to
&H^3(D_r)&\buildrel{{\pi'}^*}\over{\lra}&H^3(Y)
&\to &H^1(D_{r-1})\ ,\cr
}
$$
o\`u $c$ est le nombre de composantes connexes de $D_{r-1}$. En
particulier, le conoyau de l'injection $H^3(\iota_r,\Z)$
contient $\Z^{c-1}$: il peut
\^etre arbitrairement grand et ne peut \^etre contr\^ol\'e
seulement par une condition sur $\d(r)$.

\section{Groupes de Picard et cohomologie du faisceau
structural}\label{picard}

\ind Pour obtenir des informations sur la cohomologie de $\cO_{D_r}$
connaissant celle du faisceau $\C_{D_r}$, il est n\'ecessaire d'avoir des
informations sur les applications naturelles
$$\a^p_{D_r}:H^p(D_r,\C_{D_r})\to H^p(D_r,\cO_{D_r})\ .$$
\ind La th\'eorie de Hodge entra\^ine qu'elles sont surjectives si $D_r$ est
lisse, ce qui n'est malheureusement que rarement le cas. Par les  travaux de
du Bois, Koll\'ar et Steenbrink ([K], cor. 12.9), cela
reste vrai si
$D_r$ n'a que des singularit\'es rationnelles, ce qui n'est le cas que
dans la situation g\'en\'erique (les singularit\'es
d\'eterminantielles sont rationnelles). Nous d\'emontrons un
r\'esultat similaire avec des hypoth\`eses plus faibles;  nous dirons qu'un
sch\'ema
$X$ {\it v\'erifie la propri\'et\'e} $(P_p)$ s'il est r\'egulier en
codimension $p$ et de profondeur $>p$ en ses points ferm\'es. Un
sch\'ema normal de dimension $\ge 2$ v\'erifie $(P_1)$.

\th\label{rob}
Proposition
\enonce
Soit  $X$ une vari\'et\'e projective v\'erifiant la propri\'et\'e $(P_p)$.
L'application $\a_X^p$ est surjective  et
$$H^p(X,\cO_X)\isom \Gr^0_F\Gr^W_pH^p(X,\C)\ . $$
\endth
 
{\it D\'emonstration.} Soit $L$ un faisceau ample sur $X$. Comme $X$ est de
profondeur $>p$ en ses points ferm\'es, il existe par [G2], Exp. XII, cor.
1.4, un entier $m_0$ tel que 
$H^i(X,L^{-m})=0$ pour
$i\le p$ et
$m\ge m_0$. Le sous-sch\'ema $Y$ de $X$ d\'efini par l'annulation de $p$
sections g\'en\'erales de  
$L^{m_0}$ est r\'egulier et, si $\iota$ est l'inclusion de $Y$ dans $X$,  la
restriction
$H^p(\iota,\cO)$ est injective.

\ind Consid\'erons les complexes $\underline\Omega_X^0$ et
$\underline\Omega_Y^0$ construits par du Bois dans [dB]. On a un diagramme
commutatif
$$\matrix{ H^p(X,\C)&\buildrel{\a^p_X}\over{\lra}&H^p(X,\cO_X)
&\buildrel{\b_X^p}\over{\lra}&H^p(X,\underline\Omega_X^0)\ \ 
\cr
\noalign{\vskip 1mm}
\lda{}{}&&\lda{}{H^p(\iota,\cO)}&&\lda{}{}\cr
H^p(Y,\C)&\buildrel{\a^p_Y}\over{\lra}&H^p(Y,\cO_Y)
&\buildrel{\b_Y^p}\over{\lra}&H^p(Y,\underline\Omega_Y^0)\ ,
\cr}$$ o\`u $\b_Y^p$ est bijective car $Y$ est lisse, et $H^p(\iota,\cO)$
est injective
 comme on vient de le voir, de sorte que $\b_X^p$ est injective. Mais
$\b_X^p\circ\a^p_X$ est surjective car $X$ est propre, donc $\b^p_X$ est
surjective, donc bijective, et $\a^p_X$ est surjective.

\ind Soit $\pi:\tilde X\to X$ une d\'esingularisation de $X$.  Comme
$\iota$ se factorise \`a travers $\pi$, on a un diagramme
commutatif
$$\matrix{
&H^p(X,\C)&\buildrel{H^p(\pi,\C)}\over{\llra}&H^p(\tilde X,\C)
\cr
\noalign{\vskip 1mm}
&\lda{\a^p_X}{}&&\lda{\a^p_{\tilde X}}{}\cr
H^p(\iota,\cO):&H^p(X,\cO_X)&\buildrel{H^p(\pi,\cO)}\over{\llra}&H^p(\tilde
X,\cO_{\tilde X}) &\lra&H^p(Y,\cO_Y)\ ,
\cr}$$
qui entra\^ine que $H^p(\pi,\cO) $ est injective. Le noyau de
$H^p(\pi,\C)$ est $W_{p-1}H^p(X,\C)$ ([De], prop. 8.2.5); on en d\'eduit
$
\a^p_X(W_{p-1}H^p(X,\C))=0$ et un diagramme 
$$\matrix{
\Gr^W_pH^p(X,\C)&{\hookrightarrow}&H^p(\tilde X,\C)\ \ 
\cr
\noalign{\vskip 1mm}
\lda{\bar\a^p_X}{}&&\lda{\a^p_{\tilde X}}{}\cr
H^p(X,\cO_X)&\buildrel{H^p(\pi,\cO)}\over{\llra}&H^p(\tilde
X,\cO_{\tilde X}) \ .
\cr}$$
\ind Le noyau de
$\bar\a^p_X$ est donc
$$\Gr^W_pH^p(X,\C)\cap \Ker \a^p_{\tilde X}=\Gr^W_pH^p(X,\C)\cap F^1
H^p(\tilde X,\C)\ .$$
\ind Comme les morphismes de structures de Hodge sont stricts, le membre de
droite est $F^1\Gr^W_pH^p( X,\C)$, ce qui prouve la proposition.\cqfd

\th\label{cor}
Corollaire
\enonce
Soient  $X$ et $Y$ des vari\'et\'es projectives et $f:Y\to X$ un morphisme.
On suppose que $X$ v\'erifie la propri\'et\'e $(P_p)$.

\ind {\rm a)} Si $H^p(f,\C)$ est injective, il en est de m\^eme de
$H^p(f,\cO)$.

\ind {\rm b)} Si $Y$ v\'erifie la propri\'et\'e $(P_p)$ et que 
$H^p(f,\C)$ est bijective, il en est de m\^eme de $H^p(f,\cO)$.
\endth

{\it D\'emonstration.} Soit $\pi:\tilde Y\to Y$ une d\'esingularisation de
$Y$. La compos\'ee
$$\Gr^W_pH^p(X,\C)\lra \Gr^W_pH^p(Y,\C)\lra H^p(\tilde Y,\C) $$
est injective car le morphisme de droite l'est par [De], prop. 8.2.5.2
et le morphisme de gauche par hypoth\`ese, puisque les
morphismes de structures de Hodge mixtes sont stricts. On en d\'eduit que
l'application induite
  $$\Gr_F^0\Gr^W_pH^p(X,\C)\lra  \Gr_F^0 H^p(\tilde Y,\C)\ , $$
qui par la proposition s'identifie \`a $H^p(\pi\circ f,\cO)$, est aussi
injective, d'o\`u a). Le b) r\'esulte du fait que les
morphismes de structures de Hodge mixtes sont stricts.\cqfd

\medskip
\ind On veut maintenant montrer un r\'esultat analogue pour les groupes de
Picard. Si
$X$ est une vari\'et\'e projective, on note $\Pic(X)$ son groupe de Picard,
$\Pic^0(X)$ la composante connexe de l'\'el\'ement neutre, et $NS(X)$ le
quotient $\Pic(X)/\Pic^0(X)$; c'est un sous-groupe de $H^2(X,\Z)$, il
est  ab\'elien de type fini. Si
$X$ est normale,  $\Pic^0(X)$ est une vari\'et\'e
ab\'elienne ([G1], cor.~3.2) dont l'espace tangent \`a l'origine est
$H^1(X,\cO_X)$ et le groupe
$\Pic(X)$ est isomorphe
\`a
$\Pic^0(X)\oplus NS(X)$.

\th
Proposition
\enonce
Soient $X$ et $Y$ des  vari\'et\'es projectives irr\'eductibles, avec $X$
normale,  et
$f:Y\ra X$ un morphisme. 

\ind{\rm a)} Si  $H^1(f,\Z/\ell\Z)$ est
injective pour tout entier $\ell$ premier,  
$\Pic^0(f)$ est injective.

\ind{\rm b)} Si $H^1(f,\Z)$ est bijective et  que $H^2(f,\Z/\ell\Z)$ est
injective pour tout entier $\ell$ premier, 
$\Pic(f)$ est injective et son conoyau est sans torsion; si de plus $Y$ est
normale,
$\Pic^0(f)$ est bijective.

\ind{\rm c)} Si $X$ v\'erifie $(P_2)$, que $Y$ est normale et que
$H^1(f,\Z)$ et $H^2(f,\Z)$ sont bijectives,  
$\Pic(f)$ est bijective.
\endth
 
{\it D\'emonstration.} Pour tout entier premier
$\ell$, on a un diagramme commutatif issu des suites exactes $0\to
\mu_\ell\to{\Bbb G}_m\buildrel{(\cdot)^\ell}\over{\lra}{\Bbb G}_m\to 0$
pour $X$ et $Y$:
$$\matrix{H^0(X,\cO_X^*)&\ra&H^1(X,\Z/\ell\Z)&\ra& 
\Pic(X)&\buildrel{\times
\ell}\over{\lra} &\Pic(X)&\ra &H^2(X,\Z/\ell\Z)
\cr
\noalign{\vskip 1mm}
\lda{}{\hskip-1mm H^0(f,\cO^*)} &&\lda{}{\hskip-1mm H^1(f,\Z/\ell\Z)}
&&\lda{}{\hskip-1mm \Pic(f)} &&\lda{}{\hskip-1mm \Pic(f)}
&&\lda{}{\hskip-1mm H^2(f,\Z/\ell\Z)}\cr
H^0(Y,\cO_Y^*)&\ra&H^1(Y,\Z/\ell\Z)&\ra&
\Pic(Y)&\buildrel{\times
\ell}\over{\lra} &\Pic(Y)&\ra &H^2(Y,\Z/\ell\Z)
\cr}$$ 
\ind Les hypoth\`eses de a) entra\^inent que la multiplication
par $\ell$ est injective sur le noyau
 de  
$\Pic(f)$, donc sur celui de $\Pic^0(f)$, qui est ainsi sans torsion.
Comme c'est un sous-groupe d'une vari\'et\'e ab\'elienne, il est nul,
d'o\`u a).

\ind Sous les hypoth\`eses de b), $NS(f)$ est injective, donc aussi
$\Pic(f)$ par a). Le diagramme ci-dessus montre que la multiplication par
$\ell$ est injective sur le conoyau de $\Pic(f)$, qui est donc sans
torsion. Si $Y$ est normale, le cor.~\ref{cor} montre que l'application
tangente \`a $\Pic^0(f) $ est bijective, d'o\`u b).

\ind Sous les hypoth\`eses de c), le cor.~\ref{cor} entra\^ine que
$H^1(f,\cO)$ est bijective et $H^2(f,\cO)$ injective. Une chasse au diagramme
issu des suites exactes exponentielles $0\to 
\Z\to\cO\to\cO^*\to 0$ pour $X$ et $Y$ permet de conclure.\cqfd 

\medskip
\ind Pla\c cons-nous dans la situation du \S 1, en supposant pour simplifier
$X$ lisse. Le cor.~\ref{cor} et (\ref{rem}) entra\^inent que
$H^p(\iota_r,\cO)$ est injective pour $p\le \d(r)$. Si
$D_r$ est normal et $\d(r)\ge 2$, l'application 
$H^1(\iota_r,\cO )$  est bijective par le th.~\ref{lef} et le
cor.~\ref{cor}.  Plus g\'en\'eralement, si $\d(r-[p/2])\ge \eps(p)$ et  que
$D_r$ est non singulier en codimension $p$ et a la dimension attendue (donc
qu'il est de Cohen-Macaulay), 
$H^p(\iota_r,\cO )$ est bijective.  C'est aussi le cas pour
$p<\d(r)$ lorsque
$D_r$ est lisse et $D_{r-1}$ vide par (\ref{FL}).
 
\ind Il est tentant de conjecturer (comme dans [Ma], p.~415) que
$H^p(\iota_r,\cO )$ est bijective pour
$p<\d(r)$. Laytimi a obtenu dans [L] des r\'esultats dans ce sens en
utilisant des \theo s d'annulation; sous l'hypoth\`ese que
$D_r$ a la dimension attendue, ils entra\^inent la
conjecture lorsque $X$ est une vari\'et\'e torique ou ab\'elienne, ou
lorsque $r=\min\{e,f\}-1$, ou lorsque $e=f=r+2$ (\cf\ aussi [Ma], o\`u, sous 
des
hypoth\`eses plus restrictives sur $E$ et $F$, la
conjecture est montr\'ee avec les m\^emes m\'ethodes).

\medskip
\th\label{cor34}
Corollaire
\enonce
On se place dans la situation du \S 1, avec $X$ localement intersection
compl\`ete normale.

\ind{\rm a)} Si $\d(r)\ge 1$, l'application $\Pic^0(\iota_r)$
est injective.

\ind{\rm b)} Si $\d(r)\ge 2$, l'application $\Pic(\iota_r)$ est injective et
son conoyau est sans torsion; si de plus
$D_r$ est normale, $\Pic^0(\iota_r)$
est bijective.

\ind{\rm c)} Si $D_r$ est normale, $X$ non singuli\`ere en codimension $2$,
et $\d(r)\ge
3$,

\ind\ind $\bullet$ si $D_{r-1}$ est vide et $0<r<e$, on a $\Pic
(D_r)\isom\Pic(X)\oplus\Z [\det (K)]${\rm ;}

\ind\ind $\bullet$ si $D_{r-1}$ n'est pas vide ou si $r=0$,
l'application $\Pic(\iota_r)$
est bijective.
\endth

{\it D\'emonstration.} Les points a) et b), ainsi que le deuxi\`eme cas du
point c), d\'ecoulent de la proposition, du th.~\ref{lef} et de la
rem.~\ref{rem23}. Pour le premier cas du point c),  on applique le
raisonnement de la d\'emonstration de la proposition au diagramme commutatif
\`a lignes exactes
$$\matrix{H^1(X,\Z)&\hskip-2.7mm
\to
\hskip-2.7mm&H^1(X,\cO_X)&\hskip-2.7mm\to\hskip-2.7mm&
\Pic(X)\oplus\Z&\hskip-2.7mm\to\hskip-2.7mm
&H^2(X,\Z)\oplus\Z&\hskip-2.7mm\to\hskip-2.7mm &H^2(X,\cO_X)\cr
\lda{}{H^1(\iota_r,\Z)}&&\lda{}{H^1(\iota_r,\cO)}&&\lda{}{}&&\lda{}{}
&&\lda{}{H^2(\iota_r,\cO)}
 \cr
H^1(D_r,\Z)&\hskip-2.7mm
\to
\hskip-2.7mm&H^1(D_r,\cO_{D_r})&\hskip-2.7mm\to\hskip-2.7mm&
\Pic(D_r) &\hskip-2.7mm\to\hskip-2.7mm
&H^2(D_r,\Z) &\hskip-2.7mm\to\hskip-2.7mm
&H^2(D_r,\cO_{D_r}) \cr}
$$
\ind Ceci d\'emontre le corollaire.\cqfd 

\ind Ce r\'esultat est bien s\^ur \`a rapprocher du \theo\ de Lefschetz
de Grothendieck ([G2], Exp. XII, cor. 3.6)  qui d\'emontre le point c) dans
le cas
$e=f=1$ et
$r=0$, en toute caract\'eristique (moyennant une hypoth\`ese d'annulation
de certains groupes de cohomologie qui d\'ecoule du \theo\ de Kodaira en
caract\'eristique
$0$ et pour $X$ lisse) et sans hypoth\`ese sur les singularit\'es de
$D_0$. La m\'ethode de Grothendieck semble difficile \`a g\'en\'eraliser,
m\^eme dans le cas $E$ trivial (et $F$ ample); on aurait en effet besoin de
l'annulation des groupes de cohomologie $H^i(X,S^kF^*)$ pour $i=1,2$ et
$k>0$, alors que seuls les groupes $H^i(X,S^kF^*\otimes\det(F^*))$ sont nuls
en g\'en\'eral.
\bigskip  

\ex{Exemple} Soient $C$ une
courbe projective lisse de genre
$g$ et $d$ un entier v\'erifiant $ 3\le d\le g-1$. Si $C$ n'a pas de
$g^1_d$, on a 
$\Pic(W_d(C))\isom\Pic(JC)\oplus\Z[W_{d-1}(C)]$. Si
$C$ a un 
$g^1_d$, la
restriction
$\Pic(JC)\to \Pic(W_d(C))$ est bijective par le cor.~\ref{cor34}. Si $C$
n'est pas hyperelliptique, cela entra\^ine que pour tout point $x$ de $C$,
le diviseur de Weil
$W_{d-1}(C)+x$ de $W_d(C)$ n'est pas $\Q$-Cartier : en effet, si
$\pi:C^{(d)}\to W_d(C)$ est l'application d'Abel-Jacobi, aucun multiple 
non nul de $\pi^{-1}(W_{d-1}(C)+x)=C^{d-1}+x$ n'est dans $\pi^*\Pic(JC)$. 
En revanche, si
$C$ est hyperelliptique d'involution associ\'ee $\tau$, que $\Theta$ est un
diviseur th\^eta convenable sur
$JC$, et que $\tau(x)\ne x$,
on a
$$(\Theta+x-\tau x)\cdot W_d(C)=(g-d+1)\bigl( W_{d-1}(C)+x\bigr)\ .$$
\ind En particulier, $W_{d-1}(C)+x$ n'est pas de
Cartier dans $W_d(C)$ (puisque $\Theta$ n'est pas divisible dans
$\Pic(JC)$), mais est $\Q$-Cartier.
 
\vskip1cm
 \centerline{\bf II. Lieux de d\'eg\'en\'erescence antisym\'etriques}
\bigskip

\section{Le r\'esultat de Tu}\label{alterne}
 
\ind Soit $X$ une vari\'et\'e complexe projective irr\'eductible.
Soient 
$L$ un fibr\'e en droites et $E$ un fibr\'e vectoriel de rang $e$ sur
$X$ muni d'une forme antisym\'etrique $u:E\otimes E\ra L$, \cad\ d'une
section de $\wedge^2E^*\otimes L$, ou encore d'un morphisme antisym\'etrique
$v:E\to E^*\otimes L$. On note 
$$A_r=\{ x\in X\mid \rang (u_x)\le 2r\}
$$ et on pose
$\a(r)=
\dim(X)-{e-2r\choose 2}$. La forme antisym\'etrique $u_x$ est
de rang $\le 2r$ si et seulement s'il existe un sous-espace isotrope
de $E_x$ de dimension $e-r$. On introduit donc de nouveau
 le fibr\'e en
grassmanniennes $\pi:G=G(e-r,E)\to X$ et le lieu de ses z\'eros $Y$ de
la restriction 
$\wedge^2E^*\otimes L\to \wedge^2S^*\otimes L$.
Le morphisme $\pi$ induit par restriction un morphisme
$\pi':Y\ra A_r$ propre surjectif.
En suivant la m\'ethode de [FL], Tu d\'emontre dans [T], p. 391 (il
fait l'hypoth\`ese que $L$ est trivial, mais sa d\'emonstration marche
en g\'en\'eral) 
 que {\it  si $\wedge^2E^*\otimes L$ est ample}, $H^q(G\moins Y,\F)$
s'annule pour
$ q\ge \dim(X)+{e\choose 2}+r$. On en d\'eduit, comme dans le \S 1,
 {\it si $X$ est localement intersection compl\`ete},
$H^p(G,Y;\F)=0$ pour $p\le\a(r)$. Notons $\iota_r$ l'injection de
$A_r$ dans
$X$ et $\iota$ celle de $Y$ dans $G$. Dans le diagramme commutatif
(o\`u les groupes de cohomologie sont \`a coefficients dans $\F$)
$$\matrix{H^p(X)&\phfl{H^p(\pi)}{}&H^p(G)\ \ \cr
\pvfl{H^p(\iota_r)}{}&&\pvfl{}{H^p(\iota)}\cr
H^p(A_r)&\phfl{H^p(\pi')}{}&H^p(Y)\ ,\cr
}$$
$H^p(\pi)$ est injective; pour $p\le\a(r)$, le r\'esultat de Tu
montre que
$H^p(\iota)$ l'est aussi, et il en donc de m\^eme de
$H^p(\iota_r)$.

\medskip 

\ind On suppose $A_{r-1}$ vide. On note $c$ l'inverse dans
$H^{\bullet}(A_r,\Z )$ de la classe de Chern totale du fibr\'e vectoriel
$\Ker(v)\vert_{A_r}$ (on montrera dans la prop.~\ref{prop} que
$c_1,c_3,c_5,\ldots$ sont dans le sous-anneau
$ \iota_r^*H^\bullet(X,\Z)[c_2,c_4,\ldots]$ de
$H^\bullet(A_r,\Z)$). Puisque $\Ker(v)$ est de rang
$e-2r$, on a $\Delta_\l(c)=0$ si $\l_{e-2r+1}\ne 0$. Enfin, si $\l$
et $\mu$ sont des partitions, on notera $\l\mu$ la partition obtenue
en r\'earrangeant les parts de $\l$ et de $\mu$  en
ordre d\'ecroissant.

\th\label{ortho}
Th\'eor\`eme
\enonce
Soit $X$ une vari\'et\'e projective irr\'eductible localement intersection
compl\`ete. Soient 
$E$ un fibr\'e vectoriel et $L$ un fibr\'e en droites sur $X$, avec
$\wedge^2E^*\otimes L$ {\rm ample}, et $ E\to E^*\otimes L$ un
morphisme antisym\'etrique. Supposons $A_{r-1}=\vide$; l'application
$$\matrix{\displaystyle\bigoplus_{\l=(\l_1,\l_2,\ldots)\atop
r\ge\l_1\ge\l_2\ge\cdots \ge 0}
\hskip-5mm H^{p-4|\l|}(X,\F)&\lra&H^p(A_r,\F)\cr \sum_\l
\a_\l&\longmapsto&\sum_\l 
\Delta_{\l\l}(c)\cdot \iota_r^*\a_\l\cr}$$
 est injective pour $p\le\a(r)$,
bijective pour $p<\a(r)$. 
\endth

{\it D\'emonstration.} On a 
$$\pi'^*(\sum_\l 
\Delta_{\l\l}(c)\cdot \iota_r^*\a_\l)=\iota^*\Bigl(\sum_\l
\Delta_{\l\l}(c)\cdot\pi^*\a_\l\Bigr)\ ;$$
l'application $ \sum_\l\a_\l\mapsto\sum_\l
\Delta_{\l\l}(c)\cdot\pi^*\a_\l$, \`a valeurs dans
$H^\bullet(G,\F)$, est injective,
$\iota^*$ est injective pour $p\le \a(r)$, donc
aussi l'application du \theo . On conclut avec le lemme suivant,
o\`u les groupes de cohomologie sont \`a valeurs dans  
$\F$.\cqfd

\th
Lemme
\enonce 
On a
pour
$p<\a(r)$
$$\sum_{r\ge\l_1\ge\l_2\ge\cdots \ge 0}
h^{p-4|\l|}(X)=h^p(A_r)\ .$$
\endth

{\it D\'emonstration.} Puisque $A_{r-1}$ est vide,  $\pi':Y\ra A_r$
est le fibr\'e en grassmanniennes isotropes 
$G^0(r,\iota_r^*E/\Ker(v)) $.  Notons 
$$\cP(q)=\{ \l=(\l_1,\l_2,\ldots)\mid  r\ge\l_1\ge\l_2\ge\cdots \ge 0\
,\  |\l|=q \}$$ et 
$$\cS(q)=\{\mu =(\mu_1,\ldots,\mu_s)\mid   r\ge\mu_1>\cdots>\mu_s> 0\
,\  |\mu|=q\ , \ s\ge 0\}\ .$$
\ind On a l'\'egalit\'e (\cf\ (\ref{rel})) 
$$h^p(Y)=\sum_{q\ge 0}h^{p-2q}(A_r)\Card(\cS(q))\ .$$
\ind Proc\'edant par r\'ecurrence sur $p$, on obtient
$$\eqalign{h^p(A_r)&=h^p(Y)-\sum_{q>
0}h^{p-2q}(A_r)\Card(\cS(q))\cr
&=h^p(Y)-\sum_{q>
0,\ q'\ge
0}h^{p-2q-4q'}(X)\Card(\cS(q))\Card(\cP(q'))\ .
\cr}$$
\ind On construit une bijection entre
$\bigcup_{q+2q'=q''}\cS(q)\times\cP(q')$ et $\cP(q'')$ en associant
\`a un couple $(\mu,\l)$ la partition $\l\l\mu$, la bijection
r\'eciproque envoyant une partition de $\cP(q'')$ 
qui s'\'ecrit
$(\a_1)^{a_1}\cdots(\a_t)^{a_t}$, avec $r\ge\a_1>\cdots>\a_t>0$ sur
 la partition stricte compos\'ee des $\a_i$ pour lesquels
$a_i$ est impair et la partition
$(\a_1)^{[a_1/2]}\cdots(\a_t)^{[a_t/2]}$. On en d\'eduit
$$\eqalign{h^p(A_r)&=h^p(Y)-\sum_{q''\ge
0}h^{p-2q''}(X)\Card(\cP(q'')) 
+\sum_{q'\ge
0}h^{p-4q'}(X)\Card(\cP(q'))\cr 
&=h^p(Y)-h^p(G)+\sum_{q'\ge
0}h^{p-4q'}(X)\Card(\cP(q'))
\cr}$$
ce qui, avec le r\'esultat de Tu selon lequel $h^p(G)=h^p(Y)$, prouve
le lemme.\cqfd
 
\ex{Exemple.} {\bf Grassmanniennes de droites.} Soient $V$ un espace
vectoriel de dimension $m$ et $u:\wedge^2V^*\otimes\cO\to\cO(1)$ la
forme antisym\'etrique tautologique sur $\P\wedge^2V$; le lieu $A_1$ est
l'image du plongement de Pl\"ucker
$G=G(2,V)\hookrightarrow\P\wedge^2V$ (comparer avec \ref{segre}) et il a la
dimension attendue. Le
\theo\ donne 
$$H^p(G,\Z)\isom\bigoplus_{s=0}^{m-1}
\Delta_{(1)^{2s}}(c)\cdot H^{p-4s}(\P\wedge^2V,\Z)  $$
pour $p<\dim (G) $. Le fibr\'e vectoriel $K=\Ker(v)$ sur $G$
s'ins\`ere dans une suite exacte
$$0\to K\to \cO_G\otimes V^*\to S(1)\to 0\ ,
$$
de sorte que $c=1/c(K)=c(S(1))=1+\sigma_1+\sigma_2$, avec
$\sigma_1=c_1(\cO_G(1))$ et 
$\Delta_{(1)^{2s}}(c)\equiv\sigma_2^s$ modulo $\s_1^2$. On obtient
$H^p(G ,\Z)=0$ pour $p$ impair et
 $$H^{2q}(G,\Z)\isom\bigoplus_{0\le s\le
q/2}\sigma_1^{q-2s}\sigma_2^{4s}\Z $$
pour $2q<\dim(G)$.

\section {Un \theo\ de Lefschetz}
 
\ind On a de nouveau un \theo\ analogue au th.~\ref{lef}. Soit $m$ un entier;
on d\'efinit $\eps'(m)$ comme le reste de la division de $m$ par $4$ si
$m\ge 4$, et comme
$m+1$ si $0\le m<4$. Rappelons que $\iota_r$ d\'esigne l'injection de $A_r$
dans $X$.
  
\th\label{lefalt}
Th\'eor\`eme
\enonce
Soit $X$ une vari\'et\'e projective irr\'eductible localement intersection
compl\`ete. Soient 
$E$ un fibr\'e vectoriel et $L$ un fibr\'e en droites sur $X$, avec
$\wedge^2E^*\otimes L$ {\rm ample}, et $ E\to E^*\otimes L$ un
morphisme antisym\'etrique.
Supposons $[{m\over 4}]\le r$ et $\a(r-[{m\over 4}])\ge  \eps'(m)$;
l'application
$H^p(\iota_r,\Z)$ est bijective pour $p\le
m$.\endth

{\it D\'emonstration.} Elle suit celle du th.~\ref{lef}: on
compare de nouveau les suites spectrales pour
$\pi:G\to X$ et sa restriction
$\pi':Y\to A_r$. La fibre de
$\pi'$ au-dessus d'un point de   $A_l\moins A_{l-1} $  est la
Grassmannienne $LG$ des sous-espaces isotropes de dimension $e-r$ d'un
espace vectoriel de dimension $e$ muni d'une forme antisym\'etrique de
rang $2l$, dont la cohomologie est \'etudi\'ee dans le \S \ref{LG}. Elle
est en particulier nulle en degr\'e impair.
Les faisceaux $R^q\pi_*\Z$ et $R^q\pi'_*\Z$ sont donc nuls pour $q$ impair.
De plus, il d\'ecoule de la prop.~\ref{LG2} que la restriction
$R^{2q}\pi_*\Z\to R^{2q}\pi'_*\Z$ est surjective et que son noyau
$\cK_{2q}$ est nul sur $A_{r-[{q\over 2}]}$. La d\'emonstration proc\`ede par
r\'ecurrence sur $m$; on peut supposer $e\ge 2r+2$, auquel cas on a les
in\'egalit\'es
$$\a(t)\ge \a(t-s)+s(2s+3)\leqno{\formule}\hbox{\label{formule3}}$$
\vskip-5mm$$\eps'(t)+[{t\over 4}](2[{t\over 4}]+3)>t\
.\leqno{\formule}\hbox{\label{formule4}}$$
\ind Les \'etapes de la
d\'emonstration  sont les m\^emes que celles du th.~\ref{lef}.

{\bf Premier pas.} {\it Supposons $0\le p<m$ et $[{q\over
2}]-[{p\over 4}]\le
r$. On a} 
$$\eqalign{H^0(\cK_{2q})=0\qquad &{\rm si} \quad\a(r-[{q\over 2}])\ge
0{\rm\ ;}\cr H^{p+1}(\cK_{2q})=0\qquad &{\rm si} \quad \a(r-[{q\over
2}]-[{p\over 4}])\ge 
\eps'(p)\ .\cr}$$
  
\ind On proc\`ede comme dans la d\'emonstration du
th.~\ref{lef}. La seule chose \`a v\'erifier est que
$\a(r-[{q\over
2}]-[{p\over 4}])\ge 
\eps'(p)$ entra\^ine $\a(r-[{q\over 2}])>p $, mais cela d\'ecoule de 
(\ref{formule3}) et (\ref{formule4}).
\medskip

{\bf Deuxi\`eme pas.} {\it Supposons $[{m\over 4}]\le r$ et
$\a(r-[{m\over 4}])\ge 
\eps'(m)
$.  L'application naturelle\break
$\iota_\infty^{p,q}:{}^\pi E_\infty^{p,q}\to
{}^{\pi'}\!E_\infty^{p,q}$ est injective lorsque $p<m$ et
$p+q\le m$. D'autre part,\break
${}^{\pi'}\!E_\infty^{m,0}\isom {}^{\pi'}\!E_2^{m,0}\isom
H^m(A_r,\Z)$.}
 
\ind On proc\`ede comme dans la d\'emonstration du
th.~\ref{lef}. La seule chose \`a v\'erifier est que si $p+2q\le m$,
l'hypoth\`ese
$\a(r-[{m\over 4}])\ge  \eps'(m) $ entra\^ine  $\a(r-[{q\over 2}])\ge 0$
et, si $p+2q<m$,
$$\a(r-[{q\over 2}]-[{p\over 4}])\ge  \eps'(p)\ ;$$
si $[{q\over 2}]+[{p\over 4}]<[{m\over 4}]$, cela d\'ecoule de 
(\ref{formule3}), et si $[{q\over 2}]+[{p\over 4}]=[{m\over 4}]$, on a
$\eps'(p)\le\eps'(m)$.
\medskip

{\bf Conclusion.} On conclut comme dans la d\'emonstration du
th.~\ref{lef}, en v\'erifiant que si $p+q\le m$ et
$\a(r-[{m\over 4}])\ge\eps'(m)$, on a $\a(r)> m$, ce qui r\'esulte de 
(\ref{formule3}) et (\ref{formule4}).\cqfd

 \medskip
\ind Comme dans I, le cor.~\ref{cor}  entra\^ine  que
$H^p(\iota_r,\cO)$ est injective pour $p\le \a(r)$.
Si $\a(r-[{p\over 4}])\ge \eps'(p)$ et  que $A_r$ est non
singulier en codimension $p$ et a la dimension attendue, 
$H^p(\iota_r,\cO )$ est bijective.  C'est aussi le cas pour
$p<\a(r)$ lorsque
$D_r$ est lisse et $D_{r-1}$ vide par le th.~\ref{ortho}.
On peut conjecturer (comme dans [Ma], p.~415) que $H^p(\iota_r,\cO )$ est
bijective pour
$p<\a(r)$ (\cf\ [Ma], o\`u cette conjecture est montr\'ee sous 
des
hypoth\`eses plus restrictives).

\medskip
\th\label{cor34a}
Corollaire
\enonce
On conserve les hypoth\`eses du \theo , et on suppose en outre $X$ normale.

\ind{\rm a)} Si $\a(r)\ge 1$, l'application $\Pic^0(\iota_r)$
est injective.

\ind{\rm b)} Si $\a(r)\ge 2$, l'application $\Pic(\iota_r)$ est
injective et son conoyau est sans torsion; si de plus
$A_r$ est normale, $\Pic^0(\iota_r)$
est bijective.

\ind{\rm c)} Si $A_r$ est normale et $X$ non singuli\`ere en codimension $2$, et
 si $\a(r)\ge 3$,  l'application $\Pic(\iota_r)$
est bijective.
\endth
 
\vskip1cm
 \centerline{\bf III. Lieux de d\'eg\'en\'erescence orthogonaux}
\bigskip

\section{Un th\'eor\`eme de Bertini}

\ind Soient $X$ un sch\'ema connexe et $V$ un fibr\'e vectoriel de rang $2n$
sur
$X$ muni d'une forme quadratique
\ndeg e \`a valeurs dans un fibr\'e en droites $L$. Soient $E$ et $F$ des
sous-fibr\'es totalement isotropes maximaux de $V$. On  consid\`ere les
lieux 
$$O^r=\{\ x\in X\mid \dim(E_x\cap F_x)\ge  r\ \ {\rm et}\ \
\dim(E_x\cap F_x)\equiv r\pmod{2}\ \}\ .$$
\ind On notera que la parit\'e de $\dim(E_x\cap F_x)$ reste constante; en
particulier, soit $X=O^0$ et $O^{2r+1}=\vide$ pour tout $r$, soit $X=O^1$ et
$O^{2r}=\vide$ pour tout $r$. Le cas des lieux de d\'eg\'en\'erescence
antisym\'etriques est un cas particulier de celui-ci: si 
$v:E\to E^*\otimes L$ est un morphisme antisym\'etrique, on munit le fibr\'e
$V=E\oplus (E^*\otimes L)$  de la forme quadratique de matrice
$\pmatrix{0&1\cr 1&0\cr}$ \`a valeurs dans $L$. Les sous-fibr\'es
$E\oplus\{0\}$ et $\Im(\Id_E,v)$ sont totalement isotropes maximaux et 
$A_r=O^{e-2r}$. 
                                                                                                                                                                                                                                                                                                                                                                                                                                                                                                                                                                                                                                                                                                                                                                                                                                                 
\ind La \og codimension attendue\fg\ de $O^r$ est ${r\choose 2}$.
 On aimerait obtenir des
r\'esultats analogues \`a ceux du \S~4. Commen\c cons par
la connexit\'e (dans le cas des lieux de d\'eg\'en\'erescence
antisym\'etriques, on notera que l'hypoth\`ese \og $E^*\otimes E^*\otimes
L$ ample\fg\ de la proposition ci-dessous est plus forte que l'hypoth\`ese
\og $\wedge^2 E^*\otimes L$ ample\fg\ du r\'esultat de Tu).

\th
Proposition
\enonce Dans la situation ci-dessus, on suppose de plus que $X$ est $d$\tx
connexe, que
$E^*\otimes F^*\otimes L$ est ample, et que $X=O^k$. Pour tout
$r\ge k$ avec $r-k$ pair,   $O^r$ est 
$(d-{r\choose 2}+{k\choose 2})$\tx connexe.
 
\ind En particulier, si $X$ est irr\'eductible de dimension $>
{r\choose 2}$, le lieu $O^r$ est connexe.
 \endth

{\it D\'emonstration.} La d\'emonstration est inspir\'ee de [B].
Il suffit de traiter le cas $r=k+2$. Comme $E_x\cap
F_x$ est le noyau en $x$ de la compos\'ee $u:E \i V \ra V /F \isom
F^*\otimes L$, on a\break $O^r=D_{n-r}(u)=D_{n-r+1}(u)$, puisque la
parit\'e de
$\dim(E_x\cap F_x)$ est celle de $r$. On a $X=D_{n-r+2}(u)$; la proposition
\ref{FL2} entra\^ine alors que $O^r$ est
$$d-(r-1)^2+(r-2)^2=d-2r+3=d-{r\choose 2}+{r-2\choose 2}$$ 
connexe.\cqfd

\ind {\it Dans tout ce qui suit, on suppose} $O^{r+2}=\vide$ et on note 
$\iota_r$ l'injection de $O^r$ dans $X$; sur
$O^r$, on a une suite exacte de fibr\'es vectoriels
$$0\lra K\lra \iota_r^*E\lra
\iota_r^*F^*\otimes \iota_r^*L\lra K^*\otimes \iota_r^*L\lra 0\ ,$$
o\`u $K$ est de rang $r$. En particulier,
$\det(K)^{\otimes 2}$ est isomorphe \`a
$\iota_r^*(\det(E)\otimes\det(F)\otimes (L^*)^{\otimes r})$. Nous
allons voir que $\det(K)$ est en fait d\'ej\`a dans $\iota_r^*\Pic(X)$
en \'etudiant plus g\'en\'eralement la classe $1/c(K)=\sum c_i\in
H^\bullet(O^r,\Z)$ qui appara\^it dans (\ref{FL}).

\th\label{prop}
Proposition
\enonce
Tous les $c_i$ sont dans le sous-anneau
$ \iota_r^*A^\bullet(X)[c_2,c_4,\ldots]$ de
$A^\bullet(O^r)$.
\endth

{\it D\'emonstration.} On peut
supposer, par le principe de scindage de
[F], \S 2, que tous les fibr\'es en pr\'esence sont sommes directes
de fibr\'es en droites, et enfin, en raisonnant comme dans \loc\ p.
257, que $X$ est  un produit d'espaces projectifs. On a
$$c_m(E/K)=\sum_{i=0}^mc_i(E)c_{m-i}\quad,\quad 
c_m(F/K)=\sum_{i=0}^mc_i(F)c_{m-i}\quad{\rm et}\quad F/K\isom (E/K)^*\otimes
L\ .$$
\ind Si
$m$ est un entier impair, on a les congruences suivantes modulo
$A^\bullet(X)[c_1,\ldots,c_{m-1}]$:
$$c_m\equiv c_m(F/K)=c_m((E/K)^*\otimes L)\equiv
c_m((E/K)^*)=-c_m(E/K) \equiv -c_m\ ,
$$
\cad\ $c_m\equiv 0$ puisque $A^\bullet(X)$ est sans torsion. 
La proposition s'en d\'eduit par r\'ecurrence.\cqfd

\ex{Remarques} 1) Si $L$ est trivial, on peut montrer que pour tout entier
$m$ impair, on a dans 
$A^m(O^r)$ l'\'egalit\'e 
$\sum_{i=0}^m\iota_r^*d_i\cdot c_{m-i}=0$, o\`u les $d_i\in A^i(X)$ (avec
$d_0=1$) sont les classes introduites dans [EG].

2) Supposons $X=O^0$. Localement, on peut repr\'esenter le
morphisme associ\'e\break
$u:E \hookrightarrow  V \ra V /F \isom F^*\otimes L$ par une matrice
antisym\'etrique  dont le pfaffien d\'efinit le diviseur de Cartier
$O^2$. Sur $O^{2s}$, on a  
$\det(K) \isom  L^{\otimes s}\otimes\cO_{O^{2s}}(-O^2)$.

\bigskip

\ind Toujours sous l'hypoth\`ese $E^*\otimes F^*\otimes L$ ample, on aimerait
savoir d\'eterminer
$H^p(O^r,\Z)$ pour
$p<\dim(X)-{r\choose 2}$ par une formule analogue \`a celle du
th.~\ref{ortho} lorsque $X$ est localement intersection compl\`ete. Je ne
sais montrer qu'un cas tr\`es particulier, sous des hypoth\`eses plus fortes
(on aimerait remplacer l'hypoth\`ese de 1) par $p<\dim(X)-6 $ et celle de 2)
par $p<\dim(X)-3 $). 
 
\th\label{OO}
Proposition
\enonce
Soient $X$ un sch\'ema connexe localement intersection compl\`ete et $V$ un
fibr\'e vectoriel de rang $2n$ sur
$X$ muni d'une forme quadratique
\ndeg e \`a valeurs dans un fibr\'e en droites $L$. Soient $E$ et $F$ des
sous-fibr\'es totalement isotropes maximaux de $V$. On suppose $E^*\otimes
F^*\otimes L$ ample.

\ind 1) Supposons $X=O^0$ et $O^6=\vide$.  Soit $p$ un
entier $<\dim (X)-9$; on a 
$$H^p(O^4,\Z)\isom \iota_4^*H^p(X,\Z)\qquad si\quad
p\le 3\ ;$$ $$ H^p(O^4,\Z)\isom \iota_4^*H^p(X,\Z)\oplus
 \Z c_2\qquad si\quad p=4\ .$$

\ind 2) Supposons $X=O^1$ et $O^5=\vide$. Soit $p$ un
entier $<\dim (X)-4$; on a 
$$H^p(O^3,\Z)\isom \iota_3^*H^p(X,\Z)\qquad si\quad p\le 3\ ;$$ $$
H^p(O^3,\Z)\isom \iota_3^*H^p(X,\Z)\oplus
 \Z c_2\qquad si\quad p=4\ .$$ \endth

{\it D\'emonstration.} Pour montrer 1), faisons la construction (\ref{cons})
 avec le morphisme $E \hookrightarrow  V \ra V /F \isom
F^*\otimes L$ et $r=n-3$. Notant $S$ le fibr\'e tautologique de rang $3$
sur $G$ et $\s=1/c(S)$, on a
$$H^p(G,\Z)\isom\pi^*H^p(X,\Z)\oplus
\sigma_1\pi^*H^{p-2}(X,\Z)\qquad{\rm si}\quad p\le 3\ ;$$ 
 $$H^4(G,\Z)\isom\pi^*H^4(X,\Z)\oplus
\sigma_1\pi^*H^2(X,\Z)\oplus
 \sigma_2\pi^*H^0(X,\Z)\oplus(\sigma_1^2-\s_2)\pi^*H^0(X,\Z)\ .$$
\ind Puisque $D_{n-5}$ est vide, le sch\'ema $Y$ est un fibr\'e en
$\P^3$ au-dessus de $O^4=D_{n-4}$, d'o\`u
$$H^p(Y,\Z)\isom\pi'^*H^p(O^4,\Z)\oplus
h\pi'^*H^{p-2}(O^4,\Z)\qquad{\rm si}\quad p\le 3\ ;$$
$$H^4(Y,\Z)\isom\pi'^*H^4(O^4,\Z)\oplus
h\pi'^*H^2(O^4,\Z)\oplus h^2\pi'^*H^0(O^4,\Z)\ ,$$
o\`u $h=c_1(\cO_Y(1))$.
 On a sur $Y$ une suite exacte
$$0\lra S\vert_Y\lra \pi'^*K\lra \cO_Y(1)\lra 0$$
d'o\`u $\s\vert_Y=(1+h)\pi'^*c$ et, par la prop.~\ref{prop},
$\s_1\vert_Y\equiv h $ et $\s_2\vert_Y\equiv \pi'^*c_2$ modulo 
$\pi'^*\iota_4^*H^\bullet(X,\Z)$. Lorsque $p<\d(n-3)=\dim (X)-9$, le
r\'esultat de Fulton et Lazarsfeld du \S~1 entra\^ine $H^p(G,\Z)\isom
H^p(Y,\Z)$. On en d\'eduit 1), et 2) se montre de fa\c con analogue.\cqfd

\vskip1cm
\centerline{\bf IV. Cohomologie des grassmanniennes isotropes} 
\bigskip

\ind Soit $V$ un espace vectoriel complexe de dimension $n$ muni d'une forme
antisym\'etrique de rang $2r$. On note $LG(d,V;2r)$, ou simplement $LG$, la
grassmannienne des sous-espaces vectoriels isotropes de $V$ de dimension $d$
et
$S$ le sous-fibr\'e tautologique de rang $d$ sur  $LG $.

\section {Cas o\`u la forme antisym\'etrique est non d\'eg\'en\'er\'ee 
($n=2r$)}

\ind 
L'anneau de cohomologie de $LG$ est d\'ecrit dans [PR1], th.~1.4: si
$x_1,\ldots,x_d$ sont les racines de Chern de $S$ et
$h_0,h_1,h_2,\ldots$ les fonctions sym\'etriques homog\`enes compl\`etes
d\'efinies par la s\'erie g\'en\'eratrice
$$\sum_{j\ge 0}h_jt^j=\prod_{i=1}^d(1-tx_i)^{-1}\ ,$$
on a
$$H^{2\bullet}(LG(d,V),\Z)\isom \Z[x_1,\ldots,x_d]^{{\goth
S}_d}/(h_j(x_1^2,\ldots,x_d^2),\ j> r-d)\ .$$
\ind Cela signifie que $H^{2\bullet}(LG(d,V),\Z)$ est l'alg\`ebre engendr\'ee
par $c_1(S),\ldots,c_d(S)$, avec la relation
$${1\over c(S)c(S^*)}=0\qquad\hbox{en degr\'e}\ \ \ge 2(r-d+1)\ ,$$ qui
provient du fait que le fibr\'e $S^\bot/S$ est de rang $2(r-d)$.

\subsection On a d'autre part\label{nondeg}
$$H^{2\bullet}(G(d,V),\Z)\isom \Z[x_1,\ldots,x_d]^{{\goth
S}_d}/(h_j(x_1,\ldots,x_d),\ j> 2r-d)\ ,$$
de sorte que {\it la restriction
$$H^{2p}(G(d,V),\Z)\to H^{2p}(LG(d,V),\Z)$$ est surjective pour tout $p$ et
bijective pour} $p\le 2(r-d)+1$  (et pour tout
$p$ si $d\le 1$!), tandis que les groupes de cohomologie d'ordre impair sont
tous nuls.

\section {Cas g\'en\'eral}\label{LG}

\ind On note $k=n-2r$ la dimension du noyau $K$ de
la forme antisym\'etrique. On choisit un drapeau
$$0=V_0\i V_1\i\cdots\i V_{n-1}\i V_n=V$$
de fa\c con que   $V_k=K$ et $V_{n-j}=V_{k+j}^\bot$ pour $0\le
j\le r$, auquel est associ\'e un groupe unipotent $\Gamma$ qui agit
sur
$LG $ de fa\c con que les orbites soient ceux des sous-ensembles
suivants qui ne sont pas vides:
$$O_{\bf c}=\{\Lambda\in   LG \mid \dim (\Lambda\cap V_j)=c_j\quad{\rm
pour}\quad j=1,\ldots,n\}$$
o\`u ${\bf c}=(c_1,\ldots,c_n)$ est une suite croissante d'entiers
positifs.  On a une surjection
$$\matrix{\pi_{\bf
c}: O_{c_1,\ldots,c_n}&\longrightarrow&O_{c_{k+1}-c_k,\ldots,c_n-c_k}\cr
\Lambda&\longmapsto&(\Lambda+K)/K\cr}$$
 o\`u l'orbite de droite est dans $LG(d-c_k,V/K)$, et de nouveau une
surjection
 de $\pi_{\bf
c}^{-1}(\Lambda_0)$  sur l'orbite $O_{c_1,\ldots,c_k}$
dans la grassmannienne usuelle $G(c_k,K)$, dont la fibre en $\Lambda_1$
s'identifie \`a $\Hom(K/\Lambda_1,\Lambda_0)$. On en d\'eduit
$$\dim O_{\bf c}=\dim O_{c_{k+1}-c_k,\ldots,c_n-c_k}+\dim
O_{c_1,\ldots,c_k}+(k-c_k)(d-c_k)\ .\leqno{\formule}\hbox{\label{formuleV}}$$
\ind  Posons, pour tout entier $p\ge 0$,
$$LG_p=\bigcup_{\dim O_{\bf c}\le p}  O_{\bf c}\ .$$
\ind Chaque $\overline O_{\bf c}$ est stable par $\Gamma$; c'est donc une
r\'eunion d'orbites qui est contenue dans $O_{\bf c}\cup LG_{\dim O_{\bf
c}-1}$. Il s'ensuit que chaque $LG_p$ est ferm\'e; d'autre part, le groupe
$\Gamma$ \'etant unipotent, les orbites sont des espaces affines. On obtient
ainsi une d\'ecomposition cellulaire de $LG$ au sens de  [RX], p.~223. Les
groupes d'homologie de $LG$ d'ordre impair sont donc nuls, l'application
classe de cycles $A_\bullet(LG )\to H_{2\bullet}(LG,\Z) $ est un
isomorphisme et le groupe de Chow  de $LG$ est  ab\'elien libre engendr\'e
par les classes des adh\'erences des orbites (\loc , cor., p.~223). Cela
entra\^ine, gr\^ace \`a la formule (\ref{formuleV}),
$$A_p(LG )=\bigoplus_{c=\max(0,d-r)}^{c=\min(d,k)}\
\ \bigoplus_{p'+p''=p-(k-c)(d-c)}\bigl( A_{p'}(G(c,K))\oplus
A_{p''}(LG(d-c,V/K))\bigr)\ .
$$
\ind On peut d\'ecrire le groupe de Chow de $G=G(d,V)$ de fa\c con analogue.
Il d\'ecoule de (\ref{nondeg}) que l'application
$$\rho :A_\bullet(LG)\longrightarrow A_\bullet(G) $$
 est surjective; elle est injective en degr\'e $p$ si
$$A_{p''}(LG(d-c,V/K;2r))\to A_{p''}(G(d-c,V/K))$$ l'est pour tout 
$c$ compris entre $\max(0,d-r) $ et $\min(d,k) $ et tout $p''\le
p-(k-c)(d-c)$
\cad , par
\loc , pour
$$p\le \min_{\max(0,d-r)\le c\le \min(d-2,k)}((k-c)(d-c)+2(r-d+c)+1)$$
\cad , apr\`es un petit calcul:

\th\label{LG2} Proposition
\enonce  La restriction $H^{2p}(G(d,V),\Z)\to H^{2p}(LG(d,V;2r),\Z)$ est
surjective pour tout $p$, injective pour  $p\le 2(\dim V-d-r)+1$.
\endth

\section {Cas relatif}\label{LGrel}

\ind On a vu que la cohomologie de $LG(d,V;2r)$ est le groupe ab\'elien libre
dont les g\'en\'erateurs sont les adh\'erences des orbites non
vides de $\Gamma$. Pour
chaque suite
${\bf c}$, posons
$\lambda_j=\min\{i\mid c_i=j\}$. Par [P], pp.~173--175, l'orbite $O_{\bf c}$
n'est pas vide \ssi
$$1\le\l_1<\l_2<\cdots<\l_{c_k}\le k<\l_{c_k+1}<\cdots<\l_d\le n$$
et $\l_i-k+\l_j-k\ne 2r+1$ pour $c_k<i\le j\le d$.

\ind Dans le cas $d=r=n/2$, on obtient une bijection en associant \`a chaque
g\'en\'erateur $[\overline O_{\bf c}]$ la  suite  
$0<\l_1<\l_2<\cdots<\l_s\le r$ (o\`u $s=\max\{j\mid\l_j\le r\}$), ou encore
la partition stricte $\mu=(\mu_1,\ldots,\mu_s) $ de $r$, o\`u
$\mu_j=r+1-\l_j$; la codimension de
$\overline O_{\bf c}$ dans $LG$ est $|\mu|$ ([P], p.~177). Lorsque $s=1$, on
obtient la vari\'et\'e de Schubert sp\'eciale
$$ \{\L\in LG\mid  \L\cap V_{2r+1-\mu_1}  \ne 0\}$$ dont la classe est
$c_{\mu_1}(S^*)$. Pour chaque partition stricte
$\mu$ de $r$, il existe un polyn\^ome $\tilde Q_\mu$  \`a coefficients
entiers (d\'efini  p.~21 de [PR2]), tel que $[\overline O_{\bf c}]=\tilde
Q_\mu(c(S^*))$.

\subsection On se donne maintenant un sch\'ema $X$, un fibr\'e en droites
$L$ sur $X$ et un fibr\'e vectoriel $E$ de rang $2r$ sur
$X$ muni une forme antisym\'etrique non d\'eg\'en\'er\'ee sur $E$ \`a valeurs
dans $L$. On note $LG(r,E)$ la grassmannienne relative des sous-espaces
vectoriels isotropes maximaux des fibres de $E$ et $S$ le sous-fibr\'e
tautologique de rang $r$ sur $LG(r,E)$. L'application\label{rel}
$$\matrix{\displaystyle{ \bigoplus_{\mu=(\mu_1,\ldots,
\mu_s)\atop m\ge\mu_1>\cdots> \mu_s> 0}} H^{p-2|\mu|}(X,\Z
)&\lra&H^p(LG(r,E),\Z )\cr \sum_\mu \a_\mu&\longmapsto&\sum_\mu 
\tilde Q_\mu(c(S^*))\cdot \pi^*\a_\mu\cr}$$ est bijective ([PR2], dernier
paragraphe de la p.~22, ou [F], pp.~255--256).
 
\vskip 1cm \centerline {\pc BIBLIOGRAPHIE}\bigskip

\hangindent=1cm
[ACGH] E. Arbarello, M. Cornalba, P. Griffiths, J.  Harris,  
    Geometry of algebraic curves. I. Grundlehren 267,
Springer-Verlag, New York, 1985.

\hangindent=1cm
[B] A. Bertram, {\it An existence theorem for Prym special divisors},
Invent. Math. {\bf 90} (1987), 669--671.

\hangindent=1cm
[De] P. Deligne, {\it Th\'eorie de Hodge, III},  
Publ. Math. I.H.E.S. {\bf 44} (1974), 5--77.

\hangindent=1cm
[dB] Ph. du Bois, {\it Complexe de de Rham filtr\'e d'une vari\'et\'e
singuli\`ere}, Bull. Soc. math. France {\bf 109} (1981), 41--81.

\hangindent=1cm
[EG] D. Edidin, W. Graham, {\it Characteristic Classes and Quadric
Bundles}, Duke Math. J. {\bf 78} (1995), 277--299.

\hangindent=1cm
[E] L. Ein, {\it An Analogue of Max Noether's Theorem}, Duke Math.
J. {\bf 52} (1985), 689--706.

\hangindent=1cm
[F] W. Fulton, {\it Schubert varieties in flag bundles for the classical
groups}, Proceedings of the Hirzebruch 65 Conference on Algebraic
Geometry (Ramat Gan, 1993), 241--262, Israel Math. Conf. Proc., 9,
Bar-Ilan Univ., Ramat Gan, 1996.

\hangindent=1cm
[FL] W. Fulton, R. Lazarsfeld, {\it  On the connectedness of degeneracy
loci and special divisors}, Acta Math. {\bf 146} (1981), 
271--283.

\hangindent=1cm
[G1] A. Grothendieck, {\it Techniques de descente et
th\'eor\`emes d'existence en g\'eom\'etrie alg\'ebrique VI.
Les sch\'emas de Picard : propri\'et\'es g\'en\'erales},
S\'eminaire Bourbaki, Exp. 236, 1961/62.

\hangindent=1cm
[G2] A. Grothendieck, Cohomologie locale des faisceaux
coh\'erents et \theo s de Lefschetz locaux et globaux (SGA 2), Masson
et North Holland, Paris Amsterdam, 1968.

\hangindent=1cm [H1] H. Hamm, Lefschetz theorems for singular varieties, in
{\it Singularities, Part 1} (Arcata, Calif., 1981), 547--557, Proc. Sympos.
Pure Math. {\bf 40}, Amer. Math. Soc., Providence, R.I., 1983.

\hangindent=1cm [H2] H. Hamm, {\it Zur Homotopietyp Steinscher R\"aume}, J.
Reine Angew. Math. {\bf 338} (1983), 121--135.

\hangindent=1cm
[K] J. Koll\'ar, Shafarevich maps and automorphic forms, M. B. Porter
Lectures, Princeton University Press, Princeton, NJ, 1995.

\hangindent=1cm
[L] F. Laytimi, {\it On Degeneracy Loci}, Int.  J. Math. {\bf 7}
(1996), 745--754.

\hangindent=1cm
[M] I.	Macdonald,   {\it Symmetric Products of an Algebraic Curve},
Topology, {\bf 1} (1962), 319--343.

\hangindent=1cm
[Ma] L. Manivel, {\it Vanishing theorems for ample vector bundles},
Invent. Math. {\bf 127} (1997), 401--416.
 
\hangindent=1cm
[P] P. Pragacz,  Algebro-geometric applications of Schur $S$- and
$Q$-polynomials, in {\it Topics in invariant theory (Paris, 1989/1990)},
 Lecture Notes in Math. {\bf 1478}, Springer, Berlin, 1991, 130--191.

\hangindent=1cm
[PR1] P. Pragacz, J. Ratajski, {\it A Pieri-type theorem for
Lagrangian and odd Orthogonal Grassmannians}, J. reine angew. Math. {\bf
476} (1996), 143--189.

\hangindent=1cm
[PR2] P. Pragacz, J. Ratajski, {\it Formulas for Lagrangian and
orthogonal degeneracy loci; $\tilde Q$-polynomial approach}, Comp. Math.
{\bf 107} (1997), 11--87.

\hangindent=1cm
[RX] F. Rossell\'o Lompart, S. Xamb\'o Descamps,  Computing Chow Groups,
in {\it Algebraic Geometry (Sundance, UT, 1986)}, A. Holme, R. Speiser ed.,
Lecture Notes in Math. {\bf 1311}, Springer, Berlin-New York, 1988, 220--234.

\hangindent=1cm
[S] F. Steffen,   {\it Eine Verallgemeinerung des Krullschen
Hauptidealsatzes mit einer Anwendung auf die Brill-Noether-Theorie},
Dissertation, Bochum (1996).

\hangindent=1cm
[T] L. Tu, {\it The Connectedness of Symmetric and Skew-Symmetric
Degeneracy Loci: Even Ranks}, Trans. A.M.S. {\bf 313} (1989), 381--392.

\bye